\documentclass[11pt,reqno]{amsart}
\usepackage[a4paper,margin=1.25in]{geometry}
\usepackage{amsmath,amssymb,amsthm,mathtools,mathrsfs}
\usepackage{enumitem,booktabs,microtype}
\usepackage[dvipsnames]{xcolor}
\definecolor{darkgoldenrod}{rgb}{0.60,0.43,0.02}
\usepackage[colorlinks=true,linkcolor=darkgoldenrod,citecolor=darkgoldenrod,urlcolor=darkgoldenrod]{hyperref}
\numberwithin{equation}{section}
\newtheorem{theorem}{Theorem}[section]
\newtheorem{proposition}[theorem]{Proposition}
\newtheorem{lemma}[theorem]{Lemma}
\newtheorem{corollary}[theorem]{Corollary}
\theoremstyle{definition}
\newtheorem{definition}[theorem]{Definition}

\theoremstyle{plain}
\newtheorem{lthm}{Theorem}

\newtheorem{remark}[theorem]{Remark}

\newcommand{\F}{\mathbb F}
\newcommand{\Z}{\mathbb Z}
\newcommand{\Q}{\mathbb Q}

\newcommand{\B}{\mathcal B}
\newcommand{\GL}{\operatorname{GL}}
\newcommand{\SL}{\operatorname{SL}}

\newcommand{\Gal}{\operatorname{Gal}}
\newcommand{\Frob}{\operatorname{Frob}}
\newcommand{\Res}{\operatorname{Res}}

\newcommand{\Aut}{\operatorname{Aut}}
\newcommand{\Nm}{\operatorname N}
\newcommand{\op}{\operatorname}
\newcommand{\sep}{\mathrm{sep}}

\newcommand{\lc}{\operatorname{lc}}
\newcommand{\cM}{\mathscr M}
\newcommand{\cP}{\mathscr P}
\newcommand{\cR}{\mathscr R}
\newcommand{\cU}{\mathscr U}
\newcommand{\cV}{\mathscr V}
\newcommand{\cS}{\mathscr S}
\newcommand{\cT}{\mathscr T}

\setlist[enumerate,1]{label=\textup{(\roman*)},leftmargin=2em}
\title[Average adelic images of Drinfeld modules]{Adelic Galois representations are maximal for almost all Drinfeld modules}
\author[A. Ray]{Anwesh Ray}
\address{Chennai Mathematical Institute, H1, SIPCOT IT Park, Kelambakkam, Siruseri, Tamil Nadu 603103, India}
\email{anwesh@cmi.ac.in}
\begin{document}
\begin{abstract}
Let $r\geq 3$ be an integer, let $q$ be an odd prime
power, and set $A=\mathbb{F}_q[T]$ and $F=\mathbb{F}_q(T)$. We study the
Galois representations associated with rank $r$ Drinfeld modules
of the form
$\phi_T=T+g_1\tau+\cdots+g_r\tau^r$,
where $g_i\in A$ and $g_r\neq 0$. Ordering the coefficient
tuples $\vec{g}=(g_1,\ldots,g_r)$ by their maximum degree,
we prove that, for a density-one set of tuples, the
$\lambda$-adic Galois representation is surjective for every
prime $\lambda$ of $A$, and the adelic Galois image contains
$\mathrm{SL}_r(\widehat A)$. Moreover, for every tuple in this set,
the index of the adelic image equals the index of the
corresponding rank-one determinant representation and hence
divides $q-1$. These results extend Zywina's theorem in
rank $2$ to arbitrary ranks $r\geq 3$. As a further
consequence, we establish the existence of rank $r$
Drinfeld modules with surjective adelic Galois representations
for every $r\geq 3$, extending previous existence results
which were known for prime values of $r$.
\end{abstract}
\maketitle

\section{Introduction}
\subsection{Background and motivation}
The torsion of a Drinfeld module carries commuting actions of its
coefficient ring and of the absolute Galois group. These actions give a
compatible family of representations, one for each nonzero ideal of the
coefficient ring. These images reflect information about the Galois groups of the division
fields and the entanglement relations among them. It is expected, however not known that the adelic Galois images are maximal for many Drinfeld modules of a given rank $r$ over a global field. These Galois representations are positive characteristic analogues of the characteristic zero representations that arise from elliptic curves and abelian varieties. More broadly, Drinfeld modules are closely related to the theory of Drinfeld shtukas, which played a central role in the proofs of the global Langlands correspondence over function fields, see \cite{DrinfeldLanglands,LafforgueLanglands}.

\par In this article, we specialize to rational function fields over finite fields. Let $q$ be an odd prime power and set $A=\F_q[T]$ and $F=\F_q(T)$. A Drinfeld
module of rank $r$ over $F$ is determined by an expression
\begin{equation}\label{eq:family}
 \phi_T=T+g_1\tau+\cdots+g_r\tau^r,
\end{equation}
where $g_r\neq 0$ and $\tau$ acts on an element $x$ by $x\mapsto x^q$. For a nonzero ideal
$\mathfrak a\subset A$, the $\mathfrak a$-torsion is free of rank $r$ over
$A/\mathfrak a$. This gives a homomorphism
$\rho_{\phi,\mathfrak a}:G_F\to\GL_r(A/\mathfrak a)$, where
$G_F=\Gal(F^{\sep}/F)$. Passing to compatible inverse limits gives
\[
 \widehat\rho_\phi:G_F\longrightarrow\GL_r(\widehat A),
\]
where $\widehat{A}=\prod_{\lambda\in\Omega_A}A_\lambda$ and $\Omega_A$ is the set of nonzero prime ideals of $A$ and $A_\lambda$
is the completion of $A$ at $\lambda$.

Pink and R\"utsche proved that this adelic image is open if its geometric endomorphism ring is $A$, cf. \cite{PinkRutscheOpen}. It is natural to study the maximal possible image as one varies over all Drinfeld modules over $F$ of a given rank $r$. The analogous questions for elliptic curves over number fields
have been studied extensively from a statistical perspective. Given an elliptic curve $E$ over $\mathbb{Q}$, a prime $\ell$ is said to be exceptional if the Galois representation
\[\rho_{E, \ell}:\op{G}_{\Q}\rightarrow \op{GL}_2(\F_\ell)\]
associated to $E[\ell]$ is not surjective. Serre's open image theorem implies that a non-CM elliptic curve has only finitely many exceptional primes. A striking result of Duke
\cite{DukeExceptional} establishes that almost all elliptic
curves over $\mathbb Q$, when ordered by height, have no
exceptional primes. Jones \cite{JonesSerre} subsequently
strengthened this result by proving that almost all elliptic
curves over $\mathbb Q$ are Serre curves, meaning that the
image of their adelic Galois representation has index $2$
in $\GL_2(\widehat{\mathbb Z})$. This is the smallest
possible index for an elliptic curve over $\mathbb Q$,
owing to an unavoidable relation between its division
fields. 

Similar questions have been investigated for higher-dimensional
abelian varieties. Landesman, Swaminathan, Tao and Xu \cite{LSTXAbelian}
proved such a result for non-isotrivial rational families
of abelian varieties with big monodromy. They showed that
the adelic Galois image is as large as permitted by the
arithmetic and geometric constraints of the family for
a density-one set of specializations. Their results apply to families of
Jacobians of hyperelliptic, trigonal and plane curves,
and imply, in particular, that for every $g\geq 3$
there exist infinitely many abelian varieties over
$\mathbb Q$ whose adelic Galois image is
$\operatorname{GSp}_{2g}(\widehat{\mathbb Z})$. Zywina \cite{ZywinaAbelianFamilies} subsequently
generalized these results. For an abelian scheme over a nonempty open
subset of projective space over a number field, he proved
that the adelic Galois image of a density-one set of
specializations has uniformly bounded index in the
monodromy group of the family.

In the setting of Drinfeld modules, the author
\cite{RayAlmostAll} established an analogue of Duke's
result for the $T$-adic Galois representation. More
precisely, let $A=\F_q[T]$ and $F=\F_q(T)$, and consider
rank $r$ Drinfeld modules of the form
\[
\phi_T=T+g_1\tau+\cdots+g_r\tau^r,
\]
where $g_i\in A$ and $g_r\neq 0$. For every integer $r\geq 2$, it was shown that the
$T$-adic Galois representation
\[
\widehat\rho_{\phi,T}:G_F\longrightarrow
\GL_r(\F_q[[T]])
\]
is surjective for a density-one set of coefficient tuples
$(g_1,\ldots,g_r)$, ordered by their maximum degree.
Thus, for almost all Drinfeld modules in this family,
the Galois group of the $T^n$-division field is
$\GL_r(A/(T^n))$ for every $n\geq 1$. The proof
combines Hilbert irreducibility over function fields,
Drinfeld--Tate uniformization, and sieve methods.
However, this result concerns only the prime $(T)$
and does not address the Galois representations at
other primes of $A$, or the possible entanglement relations
between their division fields.

More recently, Zywina \cite{ZywinaAverage} obtained an
analogue of the result of Jones for
rank $2$ Drinfeld modules. For $A=\F_q[T]$, he
considered the family
\[
\phi_T=T+g_1\tau+g_2\tau^2,
\quad \text{where}\quad  g_1,g_2\in A,\quad \text{and}\quad g_2\neq 0,
\]
with the coefficients ordered by degree. He proved
that, when $q>2$, the index of the adelic Galois
image in $\GL_2(\widehat A)$ divides $q-1$ for a
density-one set of coefficient pairs. When $q=2$,
the corresponding index divides $4$. 

There has also been considerable interest in constructing
individual Drinfeld modules with surjective adelic Galois
representations. In an earlier paper, Zywina
\cite{ZywinaMaximalExample} exhibited an explicit example
of rank $2$, namely
\[
\phi_T=T+\tau-T^{q-1}\tau^2,
\]
and proved that its adelic Galois representation is
surjective whenever $q\geq 5$ is an odd prime power.
Chen \cite{ChenPrimeRank} extended this construction
to prime ranks. More precisely, for an odd prime $r$,
he considered the Drinfeld module
\[
\phi_T=T+\tau^{r-1}+T^{q-1}\tau^r
\]
and established adelic surjectivity under suitable
hypotheses on $q$, depending on $r$. These
constructions provide examples with maximal adelic
Galois image, but the restriction to prime rank in
Chen's argument leaves open the corresponding
existence question for arbitrary composite ranks.
The results of the present paper address this question
and, more generally, extend Zywina's density theorem
to every rank $r\geq 3$.

\subsection{The main results}
Throughout the paper, unless a different hypothesis is stated, $r\geq3$ is an
integer and $q=p^e$ is odd. Write $A_{\leq n}=\{a\in A:\deg a\leq n\}$, with $0$ included, and set
\[
 \B_n=A_{\leq n}^{r-1}\times(A_{\leq n}\setminus\{0\}).
\]
For a subset $S$ of the coefficient family, denote by
\[
 \underline\delta(S)=\liminf_{n\to\infty}\frac{|S\cap\B_n|}{|\B_n|}\quad\text{and}\quad\overline\delta(S)=\limsup_{n\to\infty}\frac{|S\cap\B_n|}{|\B_n|}
\]the lower and upper densities respectively. When the two values agree, we denote their common value by $\delta(S)$.

For a coefficient tuple
$\boldsymbol g=(g_1,\ldots,g_r)\in A^r$ with $g_r\neq 0$,
let $\phi^{\boldsymbol g}$ denote the associated Drinfeld
module. For each nonzero prime ideal $\lambda$ of $A$,
the action of $G_F=\Gal(F^{\mathrm{sep}}/F)$ on the
$\lambda$-adic Tate module of $\phi^{\boldsymbol g}$
gives a representation
\[
\rho_{\boldsymbol g,\lambda}:G_F
\longrightarrow \GL_r(A_\lambda),
\]
where $A_\lambda$ is the completion of $A$ at
$\lambda$. Taking the product over all nonzero primes
$\lambda$ gives the adelic Galois representation
\[
\rho_{\boldsymbol g}:G_F
\longrightarrow \GL_r(\widehat A).
\]
We denote its image by
\[
G_{\boldsymbol g}
=\rho_{\boldsymbol g}(G_F)
\subseteq \GL_r(\widehat A).
\]
For a nonzero polynomial $f\in A$, we write
$\operatorname{lc}(f)\in\F_q^\times$ for its leading coefficient.

\begin{lthm}\label{thm a}
There is a set $\cM$ of coefficient tuples of density one such that, for
$\boldsymbol g\in\cM$, the following statements hold:
\begin{enumerate}
\item $\widehat\rho_{\boldsymbol g,\lambda}(G_F)=\GL_r(A_\lambda)$ for every
$\lambda\in\Omega_A$;
\item the adelic image $G_{\boldsymbol g}$ contains $\SL_r(\widehat A)$, and
its commutator subgroup equals $\SL_r(\widehat A)$;
\item if $d=\deg g_r$ and $e_r(g_r)$ is the order of
$(-1)^{d+r-1}\lc(g_r)$ in $\F_q^\times$, then
\begin{equation}\label{eq:main-index}
 [\GL_r(\widehat A):G_{\boldsymbol g}]
 =\gcd\left(d-1,\frac{q-1}{e_r(g_r)}\right).
\end{equation}
\end{enumerate}
In particular, the adelic index divides $q-1$ on $\cM$.
\end{lthm}

For $d=0$, the greatest common divisor in \eqref{eq:main-index} is
interpreted using $\gcd(-1,b)=1$.

Let $\cP$ be the set of tuples for which
$(-1)^{\deg g_r+r-1}\lc(g_r)$ generates $\F_q^\times$, and let
$\cS_{\mathrm{ad}}$ be the set with full adelic image $\GL_r(\widehat A)$.
\begin{lthm}\label{thm b}
The set $\cP$ has density $\varphi(q-1)/(q-1)$, and
\[
 \lim_{n\to\infty}
 \frac{|\cS_{\mathrm{ad}}\cap\cP\cap\B_n|}{|\cP\cap\B_n|}=1.
\]
Consequently,
$\underline\delta(\cS_{\mathrm{ad}})\geq\varphi(q-1)/(q-1)$.
\end{lthm}

\subsection{Method of proof}

We begin by studying the Galois action at primes where
the Drinfeld module has stable reduction of rank $r-1$.
More precisely, we consider primes $\mathfrak p$ of $A$
for which
\begin{equation}\label{eq:intro-stable}
v_{\mathfrak p}(g_{r-1})=0\quad\text{and}\quad
v_{\mathfrak p}(g_r)=1.
\end{equation}
Here $v_{\mathfrak p}$ denotes the normalized valuation
at $\mathfrak p$. Drinfeld--Tate uniformization gives
a description of the module in terms of a Drinfeld module
of smaller rank and a lattice. This allows us to study
the inertia action on the division points. Under the
conditions in \eqref{eq:intro-stable}, we show that the
inertia image contains a large unipotent subgroup
(Theorem~\ref{thm:full-U}).

The next step is to study the irreducibility of the
residual representations. We show that the primes at
which reducibility can occur have bounded degree,
where the bound depends on the degrees of two
auxiliary primes of good reduction
(Theorem~\ref{thm:exceptional-bound}). These auxiliary
primes can be chosen independently of the coefficient
tuple in the families we consider. Once irreducibility
is established, the unipotent subgroups obtained from
inertia can be used to generate the special linear
group (Proposition~\ref{prop:elementary-image}).
Together with the determinant representation, this
allows us to establish surjectivity of the individual
$\lambda$-adic representations
(Proposition~\ref{prop:local-from-irreducible}).

The surjectivity of the individual $\lambda$-adic
representations does not by itself imply that the
adelic image is large. Indeed, there could be
relations between the images at distinct primes,
even when each projection is surjective. We address
this problem by studying the inertia action
simultaneously on division points of coprime orders.
This allows us to show that the images at distinct
primes are independent, in the sense that their
joint image contains the product of the corresponding
special linear groups
(Proposition~\ref{prop:pairwise}). Using the fact
that the special linear groups are perfect, we
extend this to any finite collection of primes.
It follows that the adelic image contains
$\SL_r(\widehat A)$
(Theorem~\ref{thm:adelic-independent}).
The remaining question is then to determine the
determinant image, which accounts for the possible
index in $\GL_r(\widehat A)$.

Finally, we show that the conditions needed in
these arguments hold for a density-one set of
coefficient tuples. The required conditions on
stable reduction can be expressed in terms of
congruences satisfied by the coefficients. We use
sieve theoretic arguments to show that almost
all tuples satisfy these conditions at sufficiently
many primes (Proposition~\ref{prop:stable-density}).
The remaining primes are handled using a result of Zywina from \cite{ZywinaAverage}
(see Theorem~\ref{cor:fixed-level}).
Combining these results, we obtain a density-one
set of tuples for which the preceding arguments
apply (Proposition~\ref{prop:V-density}).

\subsection*{Declarations}
\begin{description}
     \item[Conflict of interest] Not applicable; there is no conflict of interest to report.
    \item[Availability of data and materials] No data was generated or analyzed in obtaining the results in this article.
    \item[Funding] No funding was received during the preparation of this article.
    \item[Use of artificial intelligence]
The author used ChatGPT 6 Astra as an exploratory tool to study related papers and to assist in correcting mistakes in the article. It was not used to generate the mathematical ideas or the main results of the article.
\end{description}

\section{Basic notions}\label{sec:basic}
This section is devoted to setting up notation and discussing preliminaries. Let $q$ be a prime power and put $A=\F_q[T]$ and $F=\F_q(T)$. A place of $F$ is an equivalence class of valuations. These correspond to the monic irreducible
polynomials in $A$ as well as the place at $\infty$. We identify a prime ideal $\lambda$ with its monic
generator. We then set 
\[
 d_\lambda=\deg\lambda,\qquad k_\lambda=A/\lambda,
 \qquad Q_\lambda=|k_\lambda|=q^{d_\lambda}\quad\text{and}\quad A_\lambda=\varprojlim_n A/\lambda^n.
\]
The fraction field of $A_\lambda$ is $F_\lambda$, which is the completion of $F$ at $\lambda$. Its valuation is normalized
by $v_\lambda(\lambda)=1$. We use the same symbol $v_\lambda$ for the extension of this
valuation to an algebraic closure $\overline{F}$, with values in $\Q$.

The remaining place is $\infty$, and the valuation is normalized by $v_\infty(T)=-1$. The corresponding absolute value is given by $|a|_\infty=q^{\deg a}$ for $a\in A\setminus\{0\}$.
For a nonzero ideal $\mathfrak a=(a)$ with $a$ monic, set
$\deg\mathfrak a=\deg a$ and $\Nm\mathfrak a=q^{\deg a}$.

For a field $K$, let $G_K=\Gal(K^{\sep}/K)$. 
For a local field $K$, write $K^{\mathrm{un}}$ and $K^t$ for its maximal
unramified and maximal tamely ramified extensions in $K^{\sep}$. Thus
$I_K=\Gal(K^{\sep}/K^{\mathrm{un}})$ is the inertia subgroup and the wild inertia group is
$P_K=\Gal(K^{\sep}/K^t)$.

Let $q=p^e$, and let $K$ be a field containing $\F_q$. A polynomial
$P(X)\in K[X]$ is called \emph{additive} if
$P(X+Y)=P(X)+P(Y)$ in $K[X,Y]$. Equivalently, an additive polynomial
is one of the form $P(X)=\sum_i a_iX^{p^i}$, with $a_i\in K$.
An additive polynomial is called \emph{$\F_q$-linear} if
$P(cX)=cP(X)$ for every $c\in\F_q$. Equivalently, the
$\F_q$-linear polynomials are precisely those of the form
$P(X)=\sum_i a_iX^{q^i}$, where only finitely many coefficients
are nonzero.

The ring $K\{\tau\}$ of twisted polynomials consists of finite
sums $\sum_i a_i\tau^i$, with multiplication determined by
$\tau c=c^q\tau$. Consequently,
\[
(a\tau^i)(b\tau^j)=ab^{q^i}\tau^{i+j}.
\]
The correspondence $\sum_i a_i\tau^i\mapsto\sum_i a_iX^{q^i}$ identifies
multiplication with composition of additive polynomials. The coefficient $a_0$ equals the formal derivative of the corresponding
polynomial. Thus an additive polynomial with $a_0\neq0$ is separable.

An $A$-field is a field $K$ equipped with an $\F_q$-algebra homomorphism
$\gamma:A\to K$. A Drinfeld module of rank $r$ is a homomorphism
$\phi:A\to K\{\tau\}$ such that the constant coefficient of $\phi_a$ is
$\gamma(a)$ and $\deg_\tau\phi_a=r\deg a$ for $a\neq0$.
Since $A$ is a polynomial ring, it is enough to specify
$\phi_T=\gamma(T)+g_1\tau+\cdots+g_r\tau^r$ with $g_r\neq0$.
For $a=\sum_i c_iT^i$, one has that
\[\phi_a=\sum_i c_i\phi_T^{i}.\] Its
$\tau$-degree is $r\deg a$.

A Drinfeld module is said to have generic characteristic if $\gamma$ is injective. Otherwise
$\ker\gamma$ is a finite prime, called its $A$-characteristic. Over $F$ and
its completions we always use the natural embedding $A\hookrightarrow F$.
In generic characteristic every $\phi_a$ for $a\neq0$ is separable.

Let $\psi$ and $\psi$ be Drinfeld modules over the same $A$-field. A morphism $u:\phi\to\psi$ is a twisted polynomial satisfying
$u\phi_a=\psi_a u$ for all $a\in A$. A nonzero morphism is an isogeny.
The only units of $K\{\tau\}$ are $K^\times$. This is easy to see given that the $\tau$-degrees of
nonzero products add. Thus isomorphisms arise from scalars $c\in K^\times$. A Drinfeld module $\phi$ and a scalar $c$ gives rise to an isomorphic Drinfeld module $\psi$, which is defined by
\begin{equation}\label{eq:scaling}
\psi_T=c^{-1}\phi_Tc=T+\sum_{i=1}^r c^{q^i-1}g_i\tau^i.
\end{equation}

The action $a\cdot x=\phi_a(x)$ equips  $K^{\sep}$
with the structure of an $A$-module. If $a$ is prime to the $A$-characteristic,
define the $a$-torsion by \[\phi[a]=\ker\phi_a\subset K^{\sep}.\] It is an $A/(a)$-module.
The number of its elements is $q^{r\deg a}$, since $\phi_a$ is separable
of degree $q^{r\deg a}$.

\begin{lemma}\label{lem:torsion-free}
If $(a)$ is prime to the $A$-characteristic, then
$\phi[a]\simeq(A/(a))^r$ as an $A$-module.
\end{lemma}

\begin{proof}
This is \cite[Corollary~3.5.3]{Papikian}.
\end{proof}
The Galois action commutes with $A$ and gives rise to a homomorphism
\[\rho_{\phi,\mathfrak a}:G_K\to\GL_r(A/\mathfrak a).\]
The splitting field $K(\phi[\mathfrak a])$ is defined to be the fixed field of the kernel of the above homomorphism.

Let $W$ be a finite
$\F_q$-subspace of a field $L$ containing $\F_q$. We set
\[
 P_W(X)=\prod_{w\in W}(X-w)\quad\text{and}\quad
 E_W(X)=X\prod_{0\neq w\in W}(1-X/w).
\]

\begin{lemma}\label{lem:subspace-polynomial}
The polynomials $P_W$ and $E_W$ are $\F_q$-linear. Let
$W\subseteq W'$ be finite $\F_q$-subspaces. Then $E_W(W')$ is also an $\F_q$-subspace and
\begin{equation}\label{eq:finite-exponential-composition}
 E_{W'}=E_{E_W(W')}\circ E_W.
\end{equation}
\end{lemma}
\begin{proof}
We prove additivity of $P_W$ by induction on $\dim_{\F_q}W$.
For $W=0$ the assertion is immediate. Write
$W=W_0\oplus\F_qv$ with $v\notin W_0$. By induction $P_{W_0}$ is
$\F_q$-linear, and therefore
\[
 \begin{split}
 P_W(X)
 &=\prod_{c\in\F_q}P_{W_0}(X-cv)\\
 &=\prod_{c\in\F_q}\bigl(P_{W_0}(X)-cP_{W_0}(v)\bigr)\\
 &=P_{W_0}(X)^q-P_{W_0}(v)^{q-1}P_{W_0}(X).
 \end{split}
\]
This is again $\F_q$-linear. Differentiating with respect to $X$, we obtain
\[
P_W'(X)=-P_{W_0}(v)^{q-1}P_{W_0}'(X),
\]
since the derivative of $P_{W_0}(X)^q$ vanishes in characteristic $p$.
As $v\notin W_0$, we have $P_{W_0}(v)\neq 0$. Moreover,
$P_{W_0}'(X)$ is a nonzero constant by the induction hypothesis.
It follows that $P_W'(X)$ is also a nonzero constant. The product formula gives
$P_W'(0)=\prod_{0\neq w\in W}(-w)\neq0$, and
$E_W=P_W/P_W'(0)$. This proves the first assertion.

For the second, $E_W(W')$ is an $\F_q$-space of size $|W'|/|W|$,
since the kernel of $E_W$ on $W'$ is $W$. The composition on the right
of \eqref{eq:finite-exponential-composition} has degree $|W'|$,
derivative one, and root set exactly $W'$. Thus there exists a nonzero scalar $c$ such that
\[
E_{E_W(W')}\circ E_W=cE_{W'}.
\]
By definition, $E_U(0)=0$ and $E_U'(0)=1$ for every
finite $\F_q$-subspace $U$. Differentiating the preceding
identity at $X=0$ and applying the chain rule gives
\[
E_{E_W(W')}'(0)E_W'(0)=cE_{W'}'(0),
\]
and hence $c=1$.
\end{proof}
The Given a prime $\lambda$, $\phi_\lambda$ gives rise to a map $\phi[\lambda^{n+1}]\rightarrow \phi[\lambda^n]$ is
$\phi_\lambda$. It is easy to see that the map is surjective. Indeed, for $x\in\phi[\lambda^n]$, a root of
$\phi_\lambda(Y)-x$ lies in $K^{\sep}$ and is killed by
$\phi_{\lambda^{n+1}}$. Compatible bases may be chosen recursively by
lifting a basis at each stage. Consider the Tate module
\[
 T_\lambda(\phi)=\varprojlim_n\phi[\lambda^n]\simeq A_\lambda^r.
\]
The associated representation is denoted by
$\widehat\rho_{\phi,\lambda}:G_K\to\GL_r(A_\lambda)$, and its residual
representation by $\overline\rho_{\phi,\lambda}:G_K\to\GL_r(k_\lambda)$.
The adelic Tate module is \[T(\phi):=\prod_\lambda T_\lambda(\phi)\] and 
$\widehat{\rho}_\phi: G_K\to\GL_r(\widehat{A})$ be the associated big Galois representation, where $\widehat{A}=\prod_\lambda A_\lambda$. 

\section{Stable reduction and full translation inertia}\label{sec:tate}

Let $K=F_{\mathfrak p}$, let $\mathcal O_K$ be its valuation ring, and write
$v=v_{\mathfrak p}$ normalized by $v(\mathfrak p)=1$. Denote by $\mathfrak m_K$ its maximal ideal. Let $\mathbb C_K=\widehat{\overline K}$ denote the completed
algebraic closure of $K$ and extend $v$ to $\mathbb C_K$. In particular, if $L/M$ is a finite Galois extension and $M$ is algebraic over $K$, one has that
\begin{equation}\label{eq:norm-valuation}
 v\bigl(\operatorname{Norm}_{L/M}(x)\bigr)=[L:M]v(x)
\end{equation}
for $x\in L^\times$.

An integral model of $\phi$ is an expression
$\phi_T=T+\sum_i g_i\tau^i$ with $g_i\in\mathcal O_K$. It has stable
reduction of rank $s\geq1$ if $g_s$ is a unit and $g_i$ lies in the maximal
ideal for $i>s$. Its reduction then has rank $s$. We say that $\phi$ has good reduction if
$s=r$.

\begin{lemma}\label{lem:tame-stable}
Every Drinfeld module over $K$ acquires stable reduction after a tame
extension whose ramification index divides $q^s-1$ for some $1\leq s\leq r$.
\end{lemma}
\begin{proof}
Set
\[
\mu=\min_{\substack{1\leq i\leq r\\g_i\neq0}}
\frac{v(g_i)}{q^i-1},
\]
and let $s$ be the largest index for which the minimum is attained.
Write $\mu=a/b$ in lowest terms, with $b\geq1$. Since
$\mu=v(g_s)/(q^s-1)$ and $v(g_s)\in\Z$, we have
$b\mid q^s-1$. In particular, $b$ is prime to $p$.

Let $\varpi$ be a uniformizer of $K=F_{\mathfrak p}$, and choose
$u$ satisfying $u^b=\varpi$. Then $L=K(u)$ is a totally and tamely
ramified extension of degree $b$. The element $c=u^{-a}$ satisfies
$v(c)=-\mu$. The Drinfeld module $\widetilde\phi$ defined by
$\widetilde\phi_a=c^{-1}\phi_a c$ for $a\in A$ is isomorphic to
$\phi$ over $L$, via multiplication by $c$. Explicitly,
\[
 \widetilde\phi_T
 =T+\sum_{i=1}^r g_i c^{q^i-1}\tau^i.
\]
For every $i$ with $g_i\neq0$, its coefficient satisfies
\[
 v\bigl(g_i c^{q^i-1}\bigr)
 =v(g_i)-(q^i-1)\mu\geq0,
\]
with equality for $i=s$. Thus $\widetilde\phi$ has stable reduction.
\end{proof}

Let $\psi:A\to K\{\tau\}$ be a Drinfeld module over $K$.
Write $K^{\sep, \psi}$ for the additive group of $K^{\sep}$
equipped with the $A$-action $a\cdot x=\psi_a(x)$. This makes $K^{\sep, \psi}$ into an $A$-module. A $\psi$-lattice is a discrete, $G_K$-stable $A$-submodule
of $K^{\sep,\psi}$ that is free of finite rank.

For a $\psi$-lattice $\Lambda$, its associated exponential $ e_\Lambda:\mathbb C_K\longrightarrow\mathbb C_K$ is defined by
\[
 e_\Lambda(X)
 =X\prod_{0\neq\gamma\in\Lambda}(1-X/\gamma).
\]
\begin{lemma}\label{lem:lattice-exponential}
The product $e_\Lambda(z)$ converges for every $z\in\mathbb C_K$ and defines an entire
$\F_q$-linear function
$e_\Lambda:\mathbb C_K\to\mathbb C_K$ with kernel $\Lambda$.
It has an expansion
\[
 e_\Lambda(z)=z+\sum_{i\geq1}c_i z^{q^i}
\]
with coefficients $c_i\in K$
and is $G_K$-equivariant:
\[
 e_\Lambda(\sigma z)=\sigma(e_\Lambda(z))
\]
with $\sigma\in G_K$ and $z\in\mathbb C_K$.
\end{lemma}

\begin{proof}
The convergence, $\F_q$-linearity, and description of the kernel
follow from the lattice-exponential construction in
\cite[Proposition~5.1.3]{Papikian}, whose argument applies over
$\mathbb C_K$; see also \cite[Section~6.2]{Papikian}.

Choose a finite Galois extension $L/K$ containing an $A$-basis
of $\Lambda$. Since $\psi$ is defined over $K$, we have
$\Lambda\subset L$, and completeness of $L$ implies that the
coefficients of $e_\Lambda$ lie in $L$. The product formula and
$G_K$-stability of $\Lambda$ give the asserted equivariance
and show that these coefficients are fixed by $\Gal(L/K)$.
They therefore belong to $K$. For further details, see \cite[Section~6.2]{Papikian}.
\end{proof}

We use the following form of Drinfeld--Tate uniformization;
see \cite[Section~6.2, especially Theorem~6.2.11]{Papikian}.

\begin{theorem}[Drinfeld--Tate uniformization]\label{thm:uniformization}
Suppose an integral model of $\phi/K$ has stable reduction of rank $s$.
There exist a Drinfeld module $\psi/K$ of rank $s$ with good reduction
and a $\psi$-lattice $\Lambda$ of rank $r-s$ such that
\[
 e_\Lambda\circ\psi_a=\phi_a\circ e_\Lambda
\]
for $a\in A$. For every nonzero $a\in A$, the exponential induces a
$G_K$-equivariant isomorphism of $A/(a)$-modules
\begin{equation}\label{eq:uniform-torsion}
 \psi_a^{-1}(\Lambda)/\Lambda
 \xrightarrow{\ \sim\ }\phi[a],
\end{equation}
defined by $[z]\longmapsto e_\Lambda(z)$.
\end{theorem}
\begin{proof}
    For further details, see for instance \cite[Section~6.2, especially Theorem~6.2.11]{Papikian}.
\end{proof}

\begin{lemma}\label{lem:tate-sequence}
For every nonzero $a\in A$ there is a $G_K$-equivariant exact sequence
of $A/(a)$-modules
\begin{equation}\label{eq:tate-sequence}
 0\longrightarrow\psi[a]\longrightarrow\phi[a]
 \longrightarrow\Lambda/a\Lambda\longrightarrow0.
\end{equation}
\end{lemma}
\begin{proof}
We refer to \cite[Section~6.2]{Papikian} for this result.
\end{proof}

Assume $s=r-1$, $v(g_{r-1})=0$ and $v(g_r)=t>0$. In this case, $\Lambda$ has rank $1$. Choose an $A$-generator $\gamma$ of $\Lambda$. Since
$\Aut_A(\Lambda)=A^\times=\F_q^\times$, the Galois action has the form
$\sigma\gamma=\eta(\sigma)\gamma$ for a character
$\eta:G_K\to\F_q^\times$. Thus $K(\gamma)/K$ is Galois of degree
dividing $q-1$, and $\gamma\in K^t$.

\begin{lemma}\label{lem:gamma-valuation}
One has
\begin{equation}\label{eq:gamma-valuation}
 v(\gamma)=-\frac{t}{q-1}.
\end{equation}
For every nonconstant $a\in A$ and every $z$ with $\psi_a(z)=\gamma$,
\begin{equation}\label{eq:z-valuation}
 v(z)=-\frac{t}{(q-1)q^{(r-1)\deg a}}.
\end{equation}
\end{lemma}
\begin{proof}
We first show that $v(\gamma)<0$. By the integral Tate uniformization construction
\cite[Section~6.2]{Papikian}, we may choose the uniformizing
exponential so that
\[
 e_\Lambda(X)\equiv X\pmod{\mathfrak m_K}
\]
coefficientwise. Thus, we may write
\[
 e_\Lambda(X)=X+\sum_{i\geq1}c_iX^{q^i},
\]
we have $c_i\in\mathfrak m_K$ for every $i\geq1$. If $v(\gamma)\geq0$, then every nonzero higher term satisfies
\[
 v(c_i\gamma^{q^i})
 \geq 1+q^iv(\gamma)
 \geq 1+v(\gamma).
\]
Consequently, the sum of the higher terms has valuation strictly
greater than $v(\gamma)$ and cannot cancel the linear term.
This contradicts $e_\Lambda(\gamma)=0$, so $v(\gamma)<0$.

Since $\psi$ has good reduction, it follows that for every $b\in A\setminus\{0\}$,
\[
 \psi_b(X)=\sum_{i=0}^{s\deg b}u_iX^{q^i},
\]
where $u_i\in\mathcal O_K$ and $u_{s\deg b}\in\mathcal O_K^\times$.
If $v(x)<0$, then for every $i<s\deg b$ with $u_i\neq0$,
\[
 v(u_ix^{q^i})\geq q^iv(x)>q^{s\deg b}v(x).
\]
The leading term therefore has strictly smaller valuation than
every other nonzero term, and hence
\begin{equation}\label{eq:negative-evaluation}
 v(\psi_b(x))=q^{s\deg b}v(x).
\end{equation}
Every nonzero $\delta\in\Lambda$ has the form $\delta=\psi_b(\gamma)$
for some $b\neq0$. Thus
\[
 v(\delta)=q^{s\deg b}v(\gamma)\leq v(\gamma).
\]

Choose $c\in\F_q$ such that $\mathfrak p\nmid T-c$, put $b=T-c$,
and choose $z_b$ satisfying $\psi_b(z_b)=\gamma$.
Since $\psi_b$ has integral coefficients and $v(\gamma)<0$,
we must have $v(z_b)<0$. Applying \eqref{eq:negative-evaluation}
and using $\deg b=1$, we obtain
\[
 v(z_b)=\frac{v(\gamma)}{q^s}>v(\gamma).
\]
For every nonzero $\delta\in\Lambda$, it follows that
\[
 v(z_b/\delta)=v(z_b)-v(\delta)
 \geq v(z_b)-v(\gamma)>0.
\]
Each factor $1-z_b/\delta$ in the product defining $e_\Lambda(z_b)$
therefore has valuation zero. The finite partial products have
absolute value one, as does their limit. Consequently,
\[
 v(e_\Lambda(z_b))=v(z_b)=\frac{v(\gamma)}{q^s}<0.
\]
In particular, $e_\Lambda(z_b)\neq0$. We have that
\[
 \phi_b(e_\Lambda(z_b))
 =e_\Lambda(\psi_b(z_b))
 =e_\Lambda(\gamma)=0.
\]

Since $b=T-c$, we have
\[
 \phi_b(X)=bX+g_1X^q+\cdots+g_sX^{q^s}
                    +g_rX^{q^{s+1}}.
\]
Here $v(b)=v(g_s)=0$, $v(g_r)=t>0$, and all coefficients are
integral. Its Newton polygon consists of a horizontal segment
from $(1,0)$ to $(q^s,0)$, followed by the segment from
$(q^s,0)$ to $(q^{s+1},t)$. Thus every root of negative valuation
has valuation
\[
 -\frac{t}{q^{s+1}-q^s};
\]
see \cite[Section~2.5]{Papikian}. Applying this to $e_\Lambda(z_b)$
gives
\[
 \frac{v(\gamma)}{q^s}
 =-\frac{t}{q^{s+1}-q^s},
\]
and therefore $v(\gamma)=-t/(q-1)$.

Finally, let $a\in A$ be nonconstant and suppose $\psi_a(z)=\gamma$.
Integrality of the coefficients again implies $v(z)<0$, so
\eqref{eq:negative-evaluation} yields
\[
 q^{s\deg a}v(z)=v(\gamma)=-\frac{t}{q-1}.
\]
Since $s=r-1$, this proves \eqref{eq:z-valuation}.
\end{proof}

For a ring $R$, write
\[
 U(R)=\left\{\begin{pmatrix}I_{r-1}&v\\0&1\end{pmatrix}
                      :v\in R^{r-1}\right\}.
\]
\begin{theorem}\label{thm:full-U}
Assume that $\phi$ has integral coefficients, that $g_{r-1}$ is a
unit, and that $g_r$ has positive valuation $t$ not divisible by $p$. If $\mathfrak a$ is prime to
$\mathfrak p$, then, in a basis adapted to \eqref{eq:tate-sequence},
\[
 \rho_{\phi,\mathfrak a}(P_K)=U(A/\mathfrak a).
\]
In particular, for every nonzero prime $\lambda\neq\mathfrak p$, there is a compatible basis of
$T_\lambda(\phi)$ with respect to which
\[\widehat\rho_{\phi,\lambda}(P_K)=U(A_\lambda).\]
\end{theorem}
\begin{proof}
Put $s=r-1$, and let $(\psi,\Lambda)$ be the Tate uniformization
of $\phi$, where $\psi$ has good reduction of rank $s$ and
$\Lambda=A\gamma$ has rank one. Recall that the $A$-action on
$\Lambda$ is given by $\psi$, so that $a\Lambda=\psi_a(\Lambda)$.
Let $a$ be the monic generator of $\mathfrak a$. The finite-level
assertion is trivial when $a=1$, so assume that $\deg a>0$; set $R=A/(a)$ and $N=q^{s\deg a}$. Since $\psi[a]\simeq R^s$, the group $\psi[a]$ has $N$ elements.

Tate uniformization gives the $G_K$-equivariant exact sequence
\[
 0\longrightarrow\psi[a]
 \xrightarrow{\,e_\Lambda\,}\phi[a]
 \xrightarrow{\,\eta_a\,}\Lambda/a\Lambda
 \longrightarrow0;
\]
see \cite[Section~6.2, pp.~355--356]{Papikian}.
We first show that wild inertia fixes $\psi[a]$ and
$\Lambda/a\Lambda$ pointwise.
The quotient map is given by
\[
 \eta_a(e_\Lambda(x))
 =\psi_a(x)\pmod{a\Lambda},
\]
where $x\in\psi_a^{-1}(\Lambda)$. Replacing $x$ by $x+\delta$, with $\delta\in\Lambda$, changes
$\psi_a(x)$ by $\psi_a(\delta)\in a\Lambda$, so this formula
is well defined.

Good reduction of $\psi$ implies that $\psi_a$ has integral
coefficients and a unit leading coefficient. Moreover,
$\psi_a'(X)=a$ is a unit because $\mathfrak a$ is prime to
$\mathfrak p$. Its reduction therefore has degree $N$ and
$N$ distinct roots over an algebraic closure of the residue field.
Choose a finite unramified extension of $K$ whose residue field
contains these roots. Hensel's lemma lifts them to $N$ distinct
roots of $\psi_a$ in that extension. Thus
\[
 \psi[a]\subset K^{\mathrm{un}}\subset K^t,
\]
and $P_K=\Gal(K^{\sep}/K^t)$ fixes $\psi[a]$ pointwise.

The action of $G_K$ on the rank-one $A$-module $\Lambda$ gives
a continuous character
\[
 \chi_\Lambda:G_K\longrightarrow
 \Aut_A(\Lambda)=A^\times=\F_q^\times.
\]
Since $P_K$ is a pro-$p$ group and $q-1$ is prime to $p$,
the restriction of $\chi_\Lambda$ to $P_K$ is trivial.
Consequently, $\gamma\in K^t$, and $P_K$ also fixes
$\Lambda/a\Lambda$ pointwise.

Write $M=e_\Lambda(\psi[a])\subset\phi[a]$. The preceding
observations show that $P_K$ fixes $M$ pointwise and acts
trivially on $\phi[a]/M$. Hence, in any basis adapted to this
submodule, its image is contained in
\[
 U(R)=
 \left\{
 \begin{pmatrix}
 I_s&u\\
 0&1
 \end{pmatrix}
 :u\in R^s
 \right\}.
\]
It remains to prove that every column $u\in R^s$ occurs.

Choose $z\in K^{\sep}$ satisfying $\psi_a(z)=\gamma$, and set
\[
 L=K^t(z).
\]
Since $\psi_a$ is additive, its division equation has precisely
the roots
\[
 \{z+w:w\in\psi[a]\}.
\]
All these roots belong to $L$, because $\psi[a]\subset K^t$.
The polynomial $\psi_a(X)-\gamma$ is separable, its derivative
being $a\neq0$, so $L/K^t$ is finite Galois.

For $\sigma\in\Gal(L/K^t)$, define
\[
 c(\sigma)=\sigma(z)-z.
\]
Applying $\psi_a$ gives
\[
 \psi_a(c(\sigma))
 =\sigma(\gamma)-\gamma=0,
\]
so $c(\sigma)\in\psi[a]$. Since every element of $\psi[a]$
is fixed by $\Gal(L/K^t)$, we have
\[
 \begin{aligned}
 c(\sigma\tau)
 &=\sigma(\tau(z)-z)+\sigma(z)-z\\
 &=c(\tau)+c(\sigma).
 \end{aligned}
\]
Thus $c$ is a homomorphism to the additive group $\psi[a]$.
It is injective: if $c(\sigma)=0$, then $\sigma$ fixes $z$
and therefore fixes $L=K^t(z)$. In particular,
\[
 D:=[L:K^t]=|\Gal(L/K^t)|\leq|\psi[a]|=N.
\]

We now use the valuation of $z$ to prove the reverse inequality.
By \eqref{eq:z-valuation},
\[
 v(z)=-\frac{t}{(q-1)N}.
\]
Every $K$-automorphism of $K^{\sep}$ preserves the extension
of $v$, since $K$ is complete. Therefore
\[
 \begin{aligned}
 v\bigl(\operatorname{Norm}_{L/K^t}(z)\bigr)
 &=\sum_{\sigma\in\Gal(L/K^t)}v(\sigma(z))\\
 &=Dv(z)
 =-\frac{Dt}{(q-1)N}.
 \end{aligned}
\]
The norm belongs to $(K^t)^\times$. Every element of $K^t$
lies in a finite tamely ramified extension of $K$, whose
ramification index is prime to $p$. Since $v(K^\times)=\Z$,
it follows that
\[
 v\bigl((K^t)^\times\bigr)\subset
 \Z_{(p)}
 :=\left\{\frac{m}{n}:m\in\Z,\ n\geq1,\ p\nmid n\right\}.
\]
Consequently,
\[
 -\frac{Dt}{(q-1)N}\in\Z_{(p)}.
\]
Now $N$ is a power of $p$, whereas both $t$ and $q-1$
are prime to $p$. Thus all the factors of $p$ in the
denominator must be cancelled by $D$, which forces $N\mid D$.
Together with $D\leq N$, this gives $D=N$. The injective map
$c$ is therefore an isomorphism
\[
 \Gal(L/K^t)\xrightarrow{\ \sim\ }\psi[a].
\]
Since restriction gives a surjection
$P_K\to\Gal(L/K^t)$, we conclude that
\[
 \{\sigma(z)-z:\sigma\in P_K\}=\psi[a].
\]

To translate this statement into the action on $\phi[a]$, put
\[
 y=e_\Lambda(z).
\]
The functional equation gives
\[
 \phi_a(y)
 =e_\Lambda(\psi_a(z))
 =e_\Lambda(\gamma)=0,
\]
and the definition of $\eta_a$ gives
\[
 \eta_a(y)=\gamma\pmod{a\Lambda}.
\]
Thus $y$ lifts a generator of the free rank-one $R$-module
$\Lambda/a\Lambda$.

Choose an $R$-basis $w_1,\ldots,w_s$ of $\psi[a]$ and put
$v_i=e_\Lambda(w_i)$. The restriction of $e_\Lambda$ to
$\psi[a]$ is injective: its kernel is $\psi[a]\cap\Lambda$,
which is zero because $\Lambda$ is torsion-free as an
$A$-module. Hence $v_1,\ldots,v_s$ is a basis of $M$, and $v_1,\ldots,v_s,y$ is an $R$-basis of $\phi[a]$. Indeed, its images generate the
quotient, and any relation among these elements first has
zero coefficient of $y$ after applying $\eta_a$, and then has
all remaining coefficients zero by independence in $M$.

For every $\sigma\in P_K$, equivariance and additivity of
$e_\Lambda$ give
\[
 \sigma(y)-y
 =e_\Lambda(\sigma(z))-e_\Lambda(z)
 =e_\Lambda(\sigma(z)-z).
\]
As $\sigma$ varies, these differences range over all of $M$.
Since $\sigma$ fixes each $v_i$, its matrix in the chosen
basis is
\[
 \begin{pmatrix}
 I_s&u_\sigma\\
 0&1
 \end{pmatrix},
\]
where $u_\sigma$ is the coordinate column of $\sigma(y)-y$
in the basis $v_1,\ldots,v_s$. These columns range over
all of $R^s$, proving
\[
 \rho_{\phi,\mathfrak a}(P_K)=U(R).
\]

Finally, fix a nonzero prime $\lambda\neq\mathfrak p$, and
let $\ell\in A$ be its monic generator. At the levels
$a=\ell^n$, the transition maps in the torsion sequences are
$\psi_\ell$ on $\psi[\ell^{n+1}]$, $\phi_\ell$ on
$\phi[\ell^{n+1}]$, and natural reduction on
$\Lambda/\ell^{n+1}\Lambda$. These maps are compatible with
the injections and quotient maps. They therefore induce
a map
\[
 \eta:T_\lambda(\phi)\longrightarrow
 \varprojlim_n\Lambda/\ell^n\Lambda
 =\Lambda\otimes_A A_\lambda
\]
whose kernel is $T_\lambda(\psi)$.

We can exhibit a compatible lift of the generator of this
quotient. Set $z_0=\gamma$ and choose recursively
$z_{n+1}\in K^{\sep}$ such that $\psi_\ell(z_{n+1})=z_n$.
Then $\psi_{\ell^n}(z_n)=\gamma$. For $n\geq1$, put
$y_n=e_\Lambda(z_n)$. Using the identity
$\phi_a\circ e_\Lambda=e_\Lambda\circ\psi_a$ and
$y_n=e_\Lambda(z_n)$, we obtain
\[
 \begin{aligned}
 \phi_{\ell^n}(y_n)
 &=e_\Lambda\bigl(\psi_{\ell^n}(z_n)\bigr)
 =e_\Lambda(\gamma)=0,\\
 \phi_\ell(y_{n+1})
 &=e_\Lambda\bigl(\psi_\ell(z_{n+1})\bigr)
 =e_\Lambda(z_n)=y_n.
 \end{aligned}
\]
Here we have used $\psi_{\ell^n}(z_n)=\gamma$,
$\psi_\ell(z_{n+1})=z_n$, and
$\gamma\in\Lambda=\ker(e_\Lambda)$. Thus $y=(y_n)_{n\geq1}$ belongs to $T_\lambda(\phi)$, and
\[
 \eta(y)
 =(\gamma\bmod\ell^n\Lambda)_{n\geq1}
 =\gamma\otimes1.
\]
This element generates $\Lambda\otimes_A A_\lambda$, so
$\eta$ is surjective and we obtain the exact sequence
\[
 0\longrightarrow T_\lambda(\psi)
 \longrightarrow T_\lambda(\phi)
 \longrightarrow\Lambda\otimes_A A_\lambda
 \longrightarrow0.
\]
Adjoining $y$ to an $A_\lambda$-basis of the free rank-$s$
module $T_\lambda(\psi)$ gives an adapted basis of
$T_\lambda(\phi)$. Its reductions give compatible adapted
bases at every finite level.

In this basis, let
\[
 H=\widehat\rho_{\phi,\lambda}(P_K)\subset U(A_\lambda).
\]
The group $H$ is closed, being the continuous image of the
compact group $P_K$ in the Hausdorff group $U(A_\lambda)$.
The finite-level assertion shows that its reduction modulo
$\lambda^n$ is $U(A/\lambda^n)$ for every $n$. Given
$u\in U(A_\lambda)$, choose $h_n\in H$ such that
$h_n\equiv u\pmod{\lambda^n}$. Then $h_n\to u$, and closedness
of $H$ implies $u\in H$. Hence
\[
 \widehat\rho_{\phi,\lambda}(P_K)=U(A_\lambda),
\]
as required.
\end{proof}

\section{Inertia at the residual prime}\label{sec:connected}

In this section, we study the ramification of the determinant characters
associated with a non-zero $G_F$-stable subspace $ W\subsetneq\phi[\lambda]$
and its quotient $\phi[\lambda]/W$. Let $\lambda$ be a nonzero prime of $A$, let $\ell$ be its monic
generator, and set $d=\deg\ell$ and $Q=q^d$.
We work over $K=F_\lambda$, with valuation normalized by
$v(\ell)=1$. Let $\psi:A\to\mathcal O_K\{\tau\}$ be a Drinfeld module with good
reduction of rank $s$. Its $\lambda$-division polynomial can be expressed as follows:
\[
 \psi_\ell(X)=\ell X+\sum_{i=1}^{sd}b_iX^{q^i},
\]
where $b_i\in\mathcal O_K$ and $b_{sd}\in\mathcal O_K^\times$.
Let $m$ be the least index for which $b_m$ is a unit. We first
check that $d\mid m$. For every $a\in A$, consider the identity
$\psi_\ell\circ\psi_a=\psi_a\circ\psi_\ell$. Comparing the reduced coefficients of $X^{q^m}$ gives
\[
 \overline b_m\,\overline a^{\,q^m}
 =\overline a\,\overline b_m.
\]
Since $\overline b_m\neq0$ and the residues of elements of $A$
range over $\F_{q^d}$, every element of $\F_{q^d}$ is fixed
by the $q^m$-power map. Thus $d\mid m$. Write $m=hd$;
the integer $h$, which satisfies $1\leq h\leq s$, is the
height of the reduction.

Define
\[
 C=\{0\}\cup\{x\in\psi[\lambda]:v(x)>0\}.
\]
The $k_\lambda$-action on $\psi[\lambda]$ is induced by
$a\cdot x=\psi_a(x)$.

\begin{lemma}\label{lem:connected-irreducible}
The space $C$ is an $I_\lambda$-stable $k_\lambda$-subspace
of dimension $h$. The action of $I_\lambda$ on $C$ is
irreducible, and its action on $\psi[\lambda]/C$ is trivial.
\end{lemma}
\begin{proof}
The space $C$ consists of the $\lambda$-torsion points that
reduce to zero. It is therefore the subspace denoted by
$\psi[\lambda]^\circ(K^{\sep})$ in
\cite[Section~2]{PinkRutscheGroupRing}. Since $\psi$ has good reduction and $K=F_\lambda$,
Proposition~2.7 of that reference applies. Together with
the discussion preceding it, this gives
$\dim_{k_\lambda}C=h$ and shows that inertia acts trivially
on $\psi[\lambda]/C$. It also gives a $k_\lambda$-linear
isomorphism
\[
 \iota:C\xrightarrow{\sim}\F_{Q^h}
\]
under which the inertia action is multiplication by a
surjective character
\[
 \chi:I_\lambda\longrightarrow\F_{Q^h}^{\times}.
\]
Thus
\[
 \iota(\sigma x)=\chi(\sigma)\iota(x)
\]
for $\sigma\in I_\lambda$ and $x\in C$.
To prove irreducibility, let $D\subseteq C$ be a nonzero
$I_\lambda$-stable $k_\lambda$-subspace, and choose
$0\neq x\in D$. For any $0\neq y\in C$, the surjectivity
of $\chi$ gives an element $\sigma\in I_\lambda$ satisfying
\[
 \chi(\sigma)=\frac{\iota(y)}{\iota(x)}.
\]
Consequently $\iota(\sigma x)=\iota(y)$, and hence
$\sigma x=y$. Since $D$ is stable under inertia, this
implies $y\in D$. Therefore $D=C$, proving that the
inertia action on $C$ is irreducible.
\end{proof}

Set
\begin{equation}\label{eq:exponents}
 L_r=\prod_{i=1}^r(q^i-1)\quad \text{and}\quad N_r=L_r^2.
\end{equation}
For a finite-dimensional
$k_\lambda$-representation $V$ and a stable subspace $W$,
write $\det W$ and $\det(V/W)$ for the characters
\[
 \sigma\longmapsto\det(\sigma|_W)\quad \text{and}\quad
 \sigma\longmapsto\det(\sigma|_{V/W}),
\]
respectively.

\begin{lemma}\label{lem:lattice-exponent}
Let $1\leq t\leq r$, and let $M\in\GL_t(A)$. Assume that $M$ has finite
order. Then every eigenvalue of $M$, and every eigenvalue of
its reduction modulo a nonzero prime of $A$, has order
dividing $L_r$.
\end{lemma}
\begin{proof}
Choose $n\geq1$ such that $M^n=I_t$. Every eigenvalue
$\xi$ of $M$ satisfies $\xi^n=1$, so
$\xi\in\overline{\F}_q$. The coefficients of the
characteristic polynomial of $M$ are elementary symmetric
functions of its eigenvalues. They consequently lie in
\[
 A\cap\overline{\F}_q=\F_q.
\]
The minimal polynomial of $\xi$ over $\F_q$ divides this
characteristic polynomial. Its degree is therefore some
$i\leq t$, so $\xi\in\F_{q^i}^\times$ and its order
divides $q^i-1$, hence $L_r$.

Since the characteristic polynomial lies in $\F_q[X]$,
reduction modulo any nonzero prime of $A$ leaves it
unchanged. Its roots after reduction therefore satisfy
the same order bound.
\end{proof}

\begin{proposition}\label{prop:away-character}
Let $\phi/F$ have rank $r$ and let $W$ be a $G_F$-stable subspace of $V=\phi[\lambda]$.
For every finite prime $\mathfrak p\neq\lambda$, both
$(\det W)^{N_r}$ and $(\det(V/W))^{N_r}$ are trivial on
$I_{\mathfrak p}$.
\end{proposition}
\begin{proof}
Set $K=F_{\mathfrak p}$ and choose a uniformizer $\varpi$.
The construction in Lemma~\ref{lem:tame-stable} gives a
totally tamely ramified extension $K'=K(\varpi^{1/b})$ where $b\mid L_r$,
over which $\phi$ is isomorphic to an integral module
with stable reduction. This isomorphism is defined over
$K'$, so it identifies the torsion representations of
$G_{K'}$. We set $E=K'K^{\mathrm{un}}$ and
$J=\Gal(K^{\sep}/E)$. Since $K^{\mathrm{un}}$ contains
the $b$-th roots of unity, $E/K^{\mathrm{un}}$ is cyclic
of degree $b$. Thus $J$ is a normal subgroup of
$I_{\mathfrak p}$ of index $b$. In particular, we have that $\sigma^b\in J$ for all $\sigma\in I_{\mathfrak p}$.

By Tate uniformization over $K'$ there is a Drinfeld module $\psi$ with good reduction
and a $\psi$-lattice $\Lambda$ of rank at most $r-1$ for which there is the associated short exact sequence
\[
 0\longrightarrow\psi[\lambda]
 \longrightarrow V
 \longrightarrow\Lambda/\lambda\Lambda
 \longrightarrow0.
\]
Since $\mathfrak p\neq\lambda$ and $\psi$ has good
reduction, $\psi[\lambda]$ is unramified over $K'$.
The group $J\subseteq I_{K'}$ therefore fixes
$\psi[\lambda]$ pointwise.

Choose an $A$-basis $\gamma_1,\ldots,\gamma_t$ of $\Lambda$.
Since these elements belong to $K^{\sep}$, there is a finite
Galois extension $L/K'$ containing all of them.
Every element of $\Lambda$ has the form
\[
 \psi_{a_1}(\gamma_1)+\cdots+\psi_{a_t}(\gamma_t),
\]
where $a_1,\ldots,a_t\in A$.
The coefficients of each $\psi_{a_i}$ lie in $K'$, so every
element of $\Lambda$ belongs to $L$. Consequently,
$J\cap\Gal(K^{\sep}/L)$ fixes $\Lambda$ pointwise.
This subgroup has finite index in $J$, proving that the
action of $J$ on $\Lambda$ has finite image.
The intersection of these finitely many stabilizers
fixes the whole lattice, because the action commutes
with every $\psi_a$. The kernel of the lattice action
is consequently open and has finite index. Thus the action of $J$ on $\Lambda$ has finite image.

For $\sigma\in J$, its matrix on $\Lambda$ therefore
has finite order. Lemma~\ref{lem:lattice-exponent} shows
that every eigenvalue on $\Lambda/\lambda\Lambda$ has
order dividing $L_r$. The other eigenvalues on $V$
are all one, since $J$ fixes $\psi[\lambda]$.
This conclusion also holds when $\Lambda=0$.

Because $W$ is stable, a basis of $W$ extended to a
basis of $V$ gives a block upper triangular matrix for
$\sigma$. Its characteristic polynomial factors as
\[
 \det(XI_V-\sigma|_V)
 =\det(XI_W-\sigma|_W)
  \det(XI_{V/W}-\sigma|_{V/W}).
\]
Each determinant character is therefore a product of
eigenvalues whose orders divide $L_r$. Hence for $\sigma\in J$, we have that
\[
 (\det W)(\sigma)^{L_r}
 =(\det(V/W))(\sigma)^{L_r}=1.
\]

Let $\chi$ denote either of these determinant characters.
For arbitrary $\sigma\in I_{\mathfrak p}$, we have
\[
 \chi(\sigma)^{bL_r}
 =\chi(\sigma^b)^{L_r}=1.
\]
Since $b\mid L_r$, the integer $bL_r$ divides $N_r$.
Thus $\chi^{N_r}$ is trivial on $I_{\mathfrak p}$,
as claimed.
\end{proof}

\begin{proposition}\label{prop:residual-alternative}
Suppose that $\phi$ admits an integral model over
$F_\lambda$ with good reduction or stable reduction
of rank $r-1$. For every $G_F$-stable subspace
$0\neq W\neq V=\phi[\lambda]$, at least one of
$(\det W)^{N_r}$ and $(\det(V/W))^{N_r}$ is unramified
at $\lambda$.
\end{proposition}
\begin{proof}
All the local representations in this proof are taken
over $K=F_\lambda$, with $v_\lambda(\ell)=1$.
Replacing $\phi$ by the given integral model preserves
its $I_\lambda$-representation on $\lambda$-torsion.

Suppose first that $\phi$ has good reduction. Apply
Lemma~\ref{lem:connected-irreducible} to $\phi$, and
let $C\subseteq V$ be its positive-valuation torsion.
Since $C$ is irreducible under inertia, its stable
subspace $W\cap C$ is either zero or $C$.

If $W\cap C=0$, the map $W\to V/C$ is injective.
Inertia acts trivially on $V/C$, so it also acts
trivially on $W$. Thus $\det W$ is unramified.
If $C\subseteq W$, the natural map $V\to V/W$
factors through $V/C$. Hence inertia acts trivially
on $V/W$, and $\det(V/W)$ is unramified. This proves
the assertion in the good-reduction case.

Now suppose that the integral model has stable
reduction of rank $r-1$. Tate uniformization over $K$
gives a good-reduction module $\psi/K$ of rank $r-1$
and a rank-one lattice $\Lambda$, with exact sequence
\[
 0\longrightarrow\psi[\lambda]
 \xrightarrow{\,e_\Lambda\,}V
 \longrightarrow\Lambda/\lambda\Lambda
 \longrightarrow0.
\]
Here the sequence is applied at the residual prime
$\lambda$; this is permitted because $K$ has generic
$A$-characteristic. Identify $\psi[\lambda]$ with
its image in $V$, and let $C$ be the image of its
positive-valuation subspace. By
Lemma~\ref{lem:connected-irreducible}, $C$ is
irreducible under $I_\lambda$, and inertia acts
trivially on $B=\psi[\lambda]/C$.

Put $D=\Lambda/\lambda\Lambda$. The action on the
rank-one lattice takes values in
$\Aut_A(\Lambda)=\F_q^\times$. Thus the action on
the one-dimensional $k_\lambda$-space $D$ is given
by a character
\[
 \eta:I_\lambda\longrightarrow\F_q^\times
 \subseteq k_\lambda^\times,
 \qquad \eta^{q-1}=1.
\]
Quotienting the torsion sequence by $C$ gives
\[
 0\longrightarrow B\longrightarrow V/C
 \longrightarrow D\longrightarrow0.
\]
We again have $W\cap C=0$ or $C$.

If $W\cap C=0$, let $W'$ be the image of $W$ in
$V/C$, so that $W\simeq W'$ as inertia modules.
Let $D_W$ be the image of $W'$ in $D$. There is
an exact sequence
\[
 0\longrightarrow W'\cap B
 \longrightarrow W'
 \longrightarrow D_W\longrightarrow0.
\]
Inertia acts trivially on $W'\cap B$, while
$D_W$ is either zero or all of $D$. Multiplicativity
of determinants therefore gives
\[
 (\det W)|_{I_\lambda}=1
 \quad\text{or}\quad
 (\det W)|_{I_\lambda}=\eta.
\]
In either case, $(\det W)^{q-1}$ is unramified
at $\lambda$.

If $C\subseteq W$, put $Z=V/W$. The quotient map
$V/C\to Z$ is surjective. Let $B_Z$ be the image
of $B$ in $Z$. Inertia acts trivially on $B_Z$,
and $Z/B_Z$ is a quotient of $D$, hence is either
zero or isomorphic to $D$. The exact sequence
\[
 0\longrightarrow B_Z\longrightarrow Z
 \longrightarrow Z/B_Z\longrightarrow0
\]
then gives
\[
 (\det(V/W))|_{I_\lambda}=1
 \quad\text{or}\quad
 (\det(V/W))|_{I_\lambda}=\eta.
\]
Thus $(\det(V/W))^{q-1}$ is unramified at
$\lambda$ in this case. Since $q-1\mid N_r$,
both cases give the asserted alternative.
\end{proof}
\section{Frobenius polynomials and their exterior powers}\label{sec:frob}

In this section we prove estimates for Frobenius polynomials which are used to bound the degrees
of primes at which the residual representation is reducible. These bounds will be applied in
Section~\ref{sec:exceptional}.

Let $\mathfrak p$ be a prime of good reduction for $\phi$, and set
$d_{\mathfrak p}=\deg\mathfrak p$. As always, we denote by $\mathfrak p$
its monic generator, and normalize the infinite absolute value by
$|T|_\infty=q$.

\begin{theorem}\label{thm:frobenius}
There is a monic polynomial
\[
 P_{\phi,\mathfrak p}(X)=X^r+a_1X^{r-1}+\cdots+a_r\in A[X]
\]
with the following properties.
\begin{enumerate}
\item For every nonzero ideal $\mathfrak a$ prime to $\mathfrak p$,
$\rho_{\phi,\mathfrak a}$ is unramified at $\mathfrak p$, and
\[
 P_{\phi,\mathfrak p}(X)
 \equiv\det\bigl(XI-\rho_{\phi,\mathfrak a}
                         (\Frob_{\mathfrak p})\bigr)
 \pmod{\mathfrak a},
\]
where $\Frob_{\mathfrak p}$ denotes arithmetic Frobenius.
\item The constant term is $a_r=\epsilon_{\mathfrak p}\mathfrak p$
for some $\epsilon_{\mathfrak p}\in\F_q^\times$.
\item Every root $\alpha$ satisfies
$|\alpha|_\infty=q^{d_{\mathfrak p}/r}$ at every extension of
$|\cdot|_\infty$ to a field containing $\alpha$.
\end{enumerate}
\end{theorem}

\begin{proof}
The existence of $P_{\phi,\mathfrak p}\in A[X]$ and part~(i) are
recalled in \cite[Section~5 and Proposition~5.1]{ZywinaAverage}.
The constant-term formula and the coefficient bounds
\[
 |a_i|_\infty\leq q^{i d_{\mathfrak p}/r}
\]
for $1\leq i\leq r$ are given in \cite[Propositions~2.9 and~2.10]{ChenPrimeRank}.

To deduce part~(iii), fix any extension of $|\cdot|_\infty$ to a
splitting field and put $R=q^{d_{\mathfrak p}/r}$.
If a root $\alpha$ satisfied $|\alpha|_\infty>R$, then
\[
 |a_i\alpha^{r-i}|_\infty
 \leq R^i|\alpha|_\infty^{r-i}
 <|\alpha|_\infty^r
\]
for $1\leq i\leq r$. The leading term could therefore not cancel in
$P_{\phi,\mathfrak p}(\alpha)=0$. Thus all $r$ roots, counted with
multiplicity, have absolute value at most $R$. Their product has
absolute value $|a_r|_\infty=q^{d_{\mathfrak p}}=R^r$, so every
root has absolute value exactly $R$.
\end{proof}

\begin{lemma}\label{lem:root-support}
Let $L/F$ be a finite extension containing the roots of
$P_{\phi,\mathfrak p}$, and let $\beta$ be a product of a nonempty
subset of these roots, counted with multiplicity. Then $\beta$ is
integral over $A$, its valuation is zero at every finite place of $L$
not lying over $\mathfrak p$, and its valuation is positive at some
place above $\mathfrak p$.
\end{lemma}

\begin{proof}
Write the roots as $\alpha_1,\ldots,\alpha_r$. They are integral over
$A$ because $P_{\phi,\mathfrak p}$ is monic. Consequently, at every
finite place $w$ of $L$, all $v_w(\alpha_i)$ are nonnegative. Moreover,
\[
 \sum_{i=1}^r v_w(\alpha_i)
 =v_w\bigl((-1)^r a_r\bigr)=v_w(\mathfrak p).
\]
If $w$ does not lie above $\mathfrak p$, this sum is zero, so every
summand is zero. This proves integrality and the asserted vanishing
for $\beta$.

Suppose that $\beta$ is a product of $j>0$ roots. By
Theorem~\ref{thm:frobenius}, its absolute value at every place above
infinity is $q^{j d_{\mathfrak p}/r}>1$, when the absolute value
extends our normalization on $F$. Its valuation is therefore negative
at every such place. If all its finite valuations were zero, its
principal divisor would have negative degree, contradicting the fact
that principal divisors have degree zero. Some finite valuation is
therefore positive, and the preceding paragraph shows that it must
lie above $\mathfrak p$.
\end{proof}

For $1\leq j\leq r-1$, put $m_j=\binom rj$ and define
\[
 P_{\phi,\mathfrak p,j}(X)
 =\prod_{\substack{I\subseteq\{1,\ldots,r\}\\ |I|=j}}
       \left(X-\prod_{i\in I}\alpha_i\right).
\]
For any matrix $M$ over a field, the eigenvalues of $\bigwedge^j M$
are the products of $j$ eigenvalues of $M$, counted with multiplicity.
Indeed, triangularizing $M$ gives a triangular matrix on the exterior
power with precisely these diagonal entries.

Let $M_{\mathfrak p}$ be the companion matrix of
$P_{\phi,\mathfrak p}$. The preceding observation gives
\[
 P_{\phi,\mathfrak p,j}(X)
 =\det\bigl(XI-\textstyle\bigwedge^j M_{\mathfrak p}\bigr).
\]
The entries of $\bigwedge^j M_{\mathfrak p}$ are minors of
$M_{\mathfrak p}$ and hence belong to $A$. Thus
$P_{\phi,\mathfrak p,j}\in A[X]$, and this construction commutes
with reduction modulo an ideal.

For a monic polynomial $P(X)=\prod_{u=1}^m(X-\beta_u)$ and an
integer $N\geq1$, set
\[
 P^{[N]}(X)=\prod_{u=1}^m(X-\beta_u^N)
          =\Res_Y\bigl(P(Y),X-Y^N\bigr).
\]
The resultant expression shows that $P^{[N]}\in A[X]$ whenever
$P\in A[X]$, and that powering the roots also commutes with
reduction. Triangularization shows that $P^{[N]}$ is the
characteristic polynomial of $M^N$ whenever $P$ is the
characteristic polynomial of $M$. Therefore, for
$\lambda\neq\mathfrak p$,
\[
 P_{\phi,\mathfrak p,j}^{[N]}(X)
 \equiv
 \det\!\left(XI-
   \left(\textstyle\bigwedge^j
     \rho_{\phi,\lambda}(\Frob_{\mathfrak p})\right)^N\right)
 \pmod\lambda.
\]

Here is how this identity detects determinant characters. Let
$V=\phi[\lambda]$, and let $W\subset V$ be a nonzero proper
$G_F$-stable subspace. If $j=\dim_{k_\lambda}W$, then
$\bigwedge^jW$ is a stable line in $\bigwedge^jV$, with character
$\theta=\det(W)$. If instead $j=\dim_{k_\lambda}(V/W)$, then
$\bigwedge^j(V/W)$ is a one-dimensional quotient of
$\bigwedge^jV$, with character $\theta=\det(V/W)$. In either
case, a basis adapted to this line or quotient makes the Frobenius
matrix block triangular. Its characteristic polynomial has the
factor $X-\theta(\Frob_{\mathfrak p})$. Consequently,
\[
 P_{\phi,\mathfrak p,j}^{[N]}
       \bigl(\theta(\Frob_{\mathfrak p})^N\bigr)
 \equiv0\pmod\lambda.
\]

\begin{lemma}\label{lem:nonzero-resultant}
Let $\mathfrak p_1\neq\mathfrak p_2$ be primes of good reduction.
For $1\leq j_1,j_2\leq r-1$ and $N\geq1$, the polynomials
$P_{\phi,\mathfrak p_1,j_1}^{[N]}$ and
$P_{\phi,\mathfrak p_2,j_2}^{[N]}$ have no common root.
Their resultant is therefore a nonzero element of $A$.
\end{lemma}

\begin{proof}
Work in a finite extension containing all the Frobenius roots at
$\mathfrak p_1$ and $\mathfrak p_2$. A common root would give
$\beta_1^N=\beta_2^N$, where $\beta_i$ is a product of $j_i$
Frobenius roots at $\mathfrak p_i$. By
Lemma~\ref{lem:root-support}, there is a place $w$ above
$\mathfrak p_1$ for which $v_w(\beta_1)>0$. Since
$\mathfrak p_1\neq\mathfrak p_2$, the same lemma gives
$v_w(\beta_2)=0$. Taking valuations in the proposed equality
would give $Nv_w(\beta_1)=Nv_w(\beta_2)$, a contradiction.
\end{proof}

\begin{lemma}\label{lem:resultant-degree}
Let $\mathfrak p_1\neq\mathfrak p_2$ be primes of good reduction
of the same degree $d$. For $1\leq j_1,j_2\leq r-1$ and $N\geq1$,
\begin{align*}
 \deg_T\Res\bigl(P_{\phi,\mathfrak p_1,j_1}^{[N]},
                 P_{\phi,\mathfrak p_2,j_2}^{[N]}\bigr)
 &\leq\frac{Nd}{r}\max(j_1,j_2)m_{j_1}m_{j_2}\\
 &\leq Nd\max_{1\leq j\leq r-1}m_j^2.
\end{align*}
\end{lemma}

\begin{proof}
Denote the resultant by $R$; the preceding lemma gives
$R\in A\setminus\{0\}$. Let $\beta_1,\ldots,\beta_{m_{j_1}}$
be the roots of $P_{\phi,\mathfrak p_1,j_1}$, and let
$\gamma_1,\ldots,\gamma_{m_{j_2}}$ be the roots of
$P_{\phi,\mathfrak p_2,j_2}$. Fix an extension of
$|\cdot|_\infty$ to a field containing these roots.
Theorem~\ref{thm:frobenius} gives
\[
 |\beta_u^N|_\infty=q^{Nj_1d/r}\quad \text{and}\quad |\gamma_v^N|_\infty=q^{Nj_2d/r}.
\]
Since both polynomials are monic, the product formula for the
resultant and the ultrametric inequality give
\[
 |R|_\infty
 =\prod_{u=1}^{m_{j_1}}\prod_{v=1}^{m_{j_2}}
       |\beta_u^N-\gamma_v^N|_\infty
 \leq q^{(Nd/r)\max(j_1,j_2)m_{j_1}m_{j_2}}.
\]
Now $|R|_\infty=q^{\deg_T R}$, which proves the first bound.
The second follows from $\max(j_1,j_2)\leq r$ and
$m_{j_1}m_{j_2}\leq\max_j m_j^2$.
\end{proof}

\section{The determinant representation}\label{sec:det}

Once the adelic image contains $\SL_r(\widehat A)$, its index is
determined by its determinant. We identify that character with a
rank-one Drinfeld representation and apply Gekeler's image formula.
We also record the local surjectivity criterion and the elementary
count needed for Theorems~A and~B.

For $\phi_T=T+g_1\tau+\cdots+g_r\tau^r$, define
\begin{equation}\label{eq:det-module}
 \Delta=(-1)^{r-1}g_r\quad \text{and}\quad 
\psi_T=T+\Delta\tau.
\end{equation}
Thus $\psi$ is a rank-one Drinfeld module over $F$. The Weil
pairing theorem \cite[Theorem~3.7.1]{Papikian} gives compatible
$G_F$-equivariant isomorphisms
\begin{equation}\label{eq:det-isomorphism}
 \bigwedge^r_{A/\mathfrak a}\phi[\mathfrak a]
       \simeq\psi[\mathfrak a]\quad \text{and}\quad 
 \bigwedge^r_{A_\lambda}T_\lambda(\phi)
       \simeq T_\lambda(\psi).
\end{equation}
Since the action on the top exterior power is the determinant,
\[
 \det\rho_{\phi,\mathfrak a}=\rho_{\psi,\mathfrak a}\quad \text{and}\quad \det\widehat\rho_\phi=\widehat\rho_\psi.
\]

We use Gekeler's rank-one image theorem in the form stated in
\cite[Theorem~6.3]{ZywinaAverage}; see also
\cite[Theorem~3.13]{GekelerTwisted}.

\begin{theorem}[Gekeler]\label{thm:gekeler}
Let $\psi_T=T+\Delta\tau$ with $0\neq\Delta\in A$. Put
$d=\deg\Delta$, and let $e$ be the order of
$(-1)^d\lc(\Delta)$ in $\F_q^\times$. Then
\begin{align}
 [\widehat A^\times:\widehat\rho_\psi(G_F)]
 &=\gcd\left(d-1,\frac{q-1}{e}\right),
       \label{eq:gekeler-adelic}\\
 [(A/\mathfrak a)^\times:\rho_{\psi,\mathfrak a}(G_F)]
 &=\gcd\left(d-1,\frac{q-1}{e},
       \{v_{\mathfrak p}(\Delta):\mathfrak p\nmid\mathfrak a\}\right)
       \label{eq:gekeler-finite}
\end{align}
for every nonzero ideal $\mathfrak a$. In the second formula,
$\mathfrak p$ runs over the nonzero prime ideals of $A$ not
dividing $\mathfrak a$.
\end{theorem}

\begin{corollary}\label{cor:det-local}
Suppose $v_{\mathfrak p}(g_r)=1$. For every level $\mathfrak a$
prime to $\mathfrak p$,
\[
 \det\rho_{\boldsymbol g,\mathfrak a}(G_F)
       =(A/\mathfrak a)^\times.
\]
For every prime $\lambda\neq\mathfrak p$,
\[
 \det\widehat\rho_{\boldsymbol g,\lambda}(G_F)
       =A_\lambda^\times.
\]
For every coefficient tuple, without a valuation hypothesis,
the adelic determinant index divides $q-1$.
\end{corollary}

\begin{proof}
For the determinant module in \eqref{eq:det-module}, we have
$v_{\mathfrak p}(\Delta)=v_{\mathfrak p}(g_r)=1$.
If $\mathfrak p\nmid\mathfrak a$, the greatest common divisor
in \eqref{eq:gekeler-finite} therefore contains an entry equal
to one. The finite-level index is one, proving the first assertion.

If $\lambda\neq\mathfrak p$, apply this assertion to
$\mathfrak a=\lambda^n$ for every $n\geq1$.
The image of the continuous character
$\det\widehat\rho_{\boldsymbol g,\lambda}$ is compact and
hence closed in $A_\lambda^\times$. Surjectivity modulo every
$\lambda^n$ makes it dense in $A_\lambda^\times$, so it is
the whole group. Finally, \eqref{eq:gekeler-adelic} shows that
the adelic determinant index divides $(q-1)/e$, and thus $q-1$,
for every tuple.
\end{proof}

Recall that $\cP$ consists of the tuples for which
$(-1)^{\deg g_r+r-1}\lc(g_r)$ generates $\F_q^\times$.

\begin{lemma}\label{lem:primitive-density}
For every $n\geq0$,
\[
 \frac{|\cP\cap\B_n|}{|\B_n|}
       =\frac{\varphi(q-1)}{q-1}.
\]
In particular, $\delta(\cP)=\varphi(q-1)/(q-1)$.
\end{lemma}

\begin{proof}
Fix $0\leq d\leq n$. For each $c\in\F_q^\times$, there
are exactly $q^d$ polynomials $g_r$ of degree $d$ with leading
coefficient $c$. Multiplication by $(-1)^{d+r-1}$ permutes
$\F_q^\times$, which has exactly $\varphi(q-1)$ generators.
Thus, among the $(q-1)q^d$ possibilities for $g_r$ of degree
$d$, precisely $\varphi(q-1)q^d$ satisfy the defining condition
of $\cP$.

The remaining $r-1$ coefficients can be chosen freely in
$A_{\leq n}$, giving $q^{(r-1)(n+1)}$ possibilities. Hence
\begin{align*}
 |\cP\cap\B_n|
 &=q^{(r-1)(n+1)}\varphi(q-1)\sum_{d=0}^n q^d,\\
 |\B_n|
 &=q^{(r-1)(n+1)}(q-1)\sum_{d=0}^n q^d.
\end{align*}
Taking their ratio proves the assertion.
\end{proof}
\section{A uniform bounds and big image criteria}
\label{sec:exceptional}

The following theorem of Zywina will be of central importance in our density estimates. Recall that
\[
 A_{\leq n}=\{f\in\F_q[T]:\deg f\leq n\}\quad \text{and}\quad
 \B_n=A_{\leq n}^{\,r-1}\times
       \bigl(A_{\leq n}\setminus\{0\}\bigr),
\]
where $0$ is included in $A_{\leq n}$. Thus $\B_n$ consists of
coefficient tuples $\boldsymbol g=(g_1,\ldots,g_r)$ with each
coefficient of degree at most $n$ and $g_r\neq0$, defining
rank-$r$ Drinfeld modules
\[
 \phi_T=T+g_1\tau+\cdots+g_r\tau^r.
\]

\begin{theorem}\label{cor:fixed-level}
For every fixed nonzero ideal $\mathfrak a\subset A$,
\[
 \lim_{n\to\infty}
 \frac{\bigl|\{\boldsymbol g\in\B_n:
   \rho_{\boldsymbol g,\mathfrak a}(G_F)
       =\GL_r(A/\mathfrak a)\}\bigr|}{|\B_n|}=1.
\]
\end{theorem}
\begin{proof}
This is \cite[Theorem~8.1]{ZywinaAverage}, which is stated for
arbitrary rank $r\geq2$. The density there is computed in
$A_{\leq n}^r$. Since
\[
 \frac{|\B_n|}{|A_{\leq n}^r|}=1-q^{-(n+1)}\longrightarrow1,
\]
the same conclusion holds with our normalization.
\end{proof}

We combine the inertia-character results of
Section~\ref{sec:connected} with the resultant bounds of
Section~\ref{sec:frob}. The argument uses a determinant character
of an invariant subspace or quotient, so it applies to invariant
subspaces of every dimension.

\begin{remark}\label{lem:constant-character}
Let $\chi:G_F\to k_\lambda^\times$ be a continuous character
unramified at every finite prime. Then $\chi$ factors through
$\Gal(\overline{\F}_q/\F_q)$. In particular, there is
$\zeta\in k_\lambda^\times$ such that
\[
 \chi(\Frob_{\mathfrak p})=\zeta^{\deg\mathfrak p}
\]
for every finite prime $\mathfrak p$.
\end{remark}

Let $N_r$ be defined as in \eqref{eq:exponents}, and set
\begin{equation}\label{eq:C-r}
 C_r=N_r\max_{1\leq j\leq r-1}\binom rj^2.
\end{equation}

\begin{theorem}\label{thm:exceptional-bound}
Suppose $\phi/F$ has good reduction at distinct primes
$\mathfrak p_1,\mathfrak p_2$ of the same degree $d$.
Let $\lambda$ be a prime at which $\phi$ has good reduction or
stable reduction of rank $r-1$. If
$\overline\rho_{\phi,\lambda}$ is reducible, then
$\deg\lambda\leq C_rd$.
\end{theorem}
\begin{proof}
If $\lambda\in\{\mathfrak p_1,\mathfrak p_2\}$, then
$\deg\lambda=d\leq C_rd$. We may therefore assume that
$\lambda$ differs from both auxiliary primes.

Choose a nonzero proper $G_F$-stable subspace
$W\subset V=\phi[\lambda]$. By
Proposition~\ref{prop:away-character}, both
$(\det W)^{N_r}$ and $(\det(V/W))^{N_r}$ are unramified
at every finite prime other than $\lambda$.
Proposition~\ref{prop:residual-alternative} shows that at least
one of these two characters is also unramified at $\lambda$.
Choose one of the two representations $W$ and $V/W$ with
this property, denote it by $V_0$, and put
$\theta=\det V_0$ and $j=\dim_{k_\lambda}V_0$.
Thus $1\leq j\leq r-1$, and $\chi=\theta^{N_r}$ is
unramified at every finite prime.

Remark~\ref{lem:constant-character} supplies
$\zeta\in k_\lambda^\times$ with
$\theta(\Frob_{\mathfrak p})^{N_r}
=\zeta^{\deg\mathfrak p}$.
The exterior-power construction in Section~\ref{sec:frob}
therefore gives, at every good prime $\mathfrak p\neq\lambda$,
\begin{equation}\label{eq:powered-congruence}
 P_{\phi,\mathfrak p,j}^{[N_r]}
       (\zeta^{\deg\mathfrak p})\equiv0\pmod\lambda.
\end{equation}
In particular, the two monic polynomials
$P_{\phi,\mathfrak p_1,j}^{[N_r]}$ and
$P_{\phi,\mathfrak p_2,j}^{[N_r]}$ have the common root
$\zeta^d$ after reduction modulo $\lambda$.

Their resultant
\[
 R=\Res\bigl(P_{\phi,\mathfrak p_1,j}^{[N_r]},
             P_{\phi,\mathfrak p_2,j}^{[N_r]}\bigr)\in A
\]
is consequently divisible by $\lambda$. 

Lemma~\ref{lem:nonzero-resultant} implies that $R\neq0$, and
Lemma~\ref{lem:resultant-degree} asserts that
$\deg R\leq C_rd$. Since the nonzero polynomial $R$ is
divisible by $\lambda$, we obtain
$\deg\lambda\leq\deg R\leq C_rd$.
\end{proof}

Suppose $\mathfrak p\neq\lambda$ is a prime at which $\phi$ has
integral coefficients and
\[
 v_{\mathfrak p}(g_{r-1})=0\quad \text{and}\quad
 v_{\mathfrak p}(g_r)=1.
\]
By Theorem~\ref{thm:full-U}, there is a basis of
$T_\lambda(\phi)$ adapted to the uniformization sequence in which
\[
 \widehat\rho_{\phi,\lambda}(P_{F_{\mathfrak p}})
 =
 U(A_\lambda)
 =
 \left\{
 \begin{pmatrix}
 I_{r-1}&u\\
 0&1
 \end{pmatrix}
 :u\in A_\lambda^{r-1}
 \right\}.
\]
Thus the global image contains every matrix in this translation
group. We show that residual irreducibility then forces it to
contain $\SL_r(A_\lambda)$. 

For $i\neq j$ and $a$ in a commutative ring $R$, write
$E_{ij}(a)=I+ae_{ij}$, where $e_{ij}$ is the corresponding
matrix unit. We use the notation $[x,y]=xyx^{-1}y^{-1}$.

\begin{lemma}\label{lem:elementary-generation}
For every $r\geq2$, there is an integer $B(r)$ such that,
for every commutative local ring $R$, every element of
$\SL_r(R)$ is a product of at most $B(r)$ elementary matrices.
\end{lemma}
\begin{proof}
By \cite[Theorem~1 and Section~1]{VavilovSmolenskySury},
every element of $\SL_r(R)$ is a product of four matrices,
each upper or lower triangular with all diagonal entries
equal to one. Each such factor is a product of at most
$r(r-1)/2$ elementary matrices. This gives a bound depending
only on $r$.
\end{proof}
\begin{lemma}\label{lem:perfect}
For a commutative local ring $R$ and $r\geq3$, the group
$\SL_r(R)$ is perfect.
Moreover, $\SL_r(\widehat A)$ is perfect.
\end{lemma}
\begin{proof}
Since $r\geq3$, for every pair $i\neq j$ we can choose an index
$k$ distinct from both. We find that
\begin{equation}\label{eq:steinberg}
 [E_{ik}(a),E_{kj}(b)]=E_{ij}(ab).
\end{equation}
Taking $b=1$ shows that every elementary matrix belongs to
the commutator subgroup of $\SL_r(R)$. By
Lemma~\ref{lem:elementary-generation}, elementary matrices
generate $\SL_r(R)$, so this group is perfect.

We now prove the assertion for
$\widehat A=\prod_\lambda A_\lambda$.
Let $M=(M_\lambda)_\lambda\in\SL_r(\widehat A)$, and put
$B=B(r)$. Each $M_\lambda$ is a product of at most $B$
elementary matrices. List the pairs $(i,j)$ with $i\neq j$ as
\[
 (i_1,j_1),\ldots,(i_m,j_m), \quad \text{where}\quad m=r(r-1),
\]
and write
\[
 \mathcal E(x_1,\ldots,x_m)
 =E_{i_1j_1}(x_1)\cdots E_{i_mj_m}(x_m).
\]
Every elementary matrix can be written in this form by taking
its coefficient in the corresponding position and setting all
other parameters equal to zero. The identity matrix is obtained
by setting every parameter equal to zero.

Consequently, after adding identity factors when necessary,
we can write each local component as
\[
 M_\lambda
 =\prod_{t=1}^{B}
   \mathcal E(a_{\lambda,t,1},\ldots,a_{\lambda,t,m}),
 \qquad a_{\lambda,t,s}\in A_\lambda.
\]
For each fixed $t$ and $s$, form the element
\[
 a_{t,s}=(a_{\lambda,t,s})_\lambda\in\widehat A.
\]
We then have
\[
 M=\prod_{t=1}^{B}\mathcal E(a_{t,1},\ldots,a_{t,m}).
\]
Indeed, projection to the $\lambda$-coordinate gives exactly
the displayed factorization of $M_\lambda$.

Thus $M$ is a product of at most $Bm$ elementary matrices
over $\widehat A$. The bound is independent of $\lambda$,
so this is a finite product even though $\widehat A$ has
infinitely many local factors. Applying
\eqref{eq:steinberg} with $b=1$ over $\widehat A$ expresses
each factor as a commutator in $\SL_r(\widehat A)$.
Every element of $\SL_r(\widehat A)$ is therefore a finite
product of commutators, proving perfectness.
\end{proof}

\par Before stating the next result, we introduce some further notation. Let $R$ be a commutative local ring and put $V=R^r$.
An $R$-linear functional
$f\in V^*=\operatorname{Hom}_R(V,R)$ is called primitive
if it is surjective.

For a primitive functional $f$ and a vector $v\in\ker f$,
write $vf$ for the $R$-linear endomorphism of $V$ defined by
\[
 (vf)(x)=f(x)v.
\]
Thus $I+vf$ sends $x$ to $x+f(x)v$, where $I$ denotes
the identity map on $V$. Define
\[
 U_f=\{I+vf:v\in\ker f\}.
\]
For $v,w\in\ker f$, the relation $f(w)=0$ gives
$(vf)(wf)=0$. Consequently,
\[
 (I+vf)(I+wf)=I+(v+w)f\quad\text{and}\quad
 (I+vf)^{-1}=I-vf.
\]
These identities show that $U_f$ is a subgroup of $\GL_r(R)$.
Since $f$ is surjective, the map \[v\mapsto I+vf\] is injective.
This identifies the additive group $\ker f$ with $U_f$.

We describe this group in an adapted basis. Let
$e_1,\ldots,e_r$ be the standard basis of $V$, and write
\[
 f(x_1,\ldots,x_r)=a_1x_1+\cdots+a_rx_r.
\]
Surjectivity means that the coefficients $a_i$ generate $R$.
Since $R$ is local, at least one $a_i$ is a unit.
After reordering the basis, assume $a_r\in R^\times$, and set $b_i=e_i-a_i a_r^{-1}e_r$ for $1\leq i<r$ and set $b_r=a_r^{-1}e_r$. The change-of-basis matrix has determinant $a_r^{-1}$,
so $b_1,\ldots,b_r$ form a basis of $V$. Moreover, $f(b_i)=0$ for $1\leq i<r$ and $f(b_r)=1$. Thus $b_1,\ldots,b_{r-1}$ form a basis of $\ker f$.
If $v=\sum_{i=1}^{r-1}u_i b_i$, then $I+vf$ fixes
$b_1,\ldots,b_{r-1}$ and sends $b_r$ to $b_r+v$.
Its matrix in this basis is therefore
\[
 \begin{pmatrix}
 I_{r-1}&u\\
 0&1
 \end{pmatrix}
\]where $u=(u_1,\ldots,u_{r-1})^{\mathsf T}$. Consequently, in this basis $U_f$ is exactly the group
$U(R)$ appearing in Theorem~\ref{thm:full-U}.
In particular, $U_f\subseteq\SL_r(R)$.

Finally, let $g\in\GL_r(R)$. For every $x\in V$,
\[
 g(I+vf)g^{-1}(x)
 =x+f(g^{-1}x)\,gv.
\]
Hence
\[
 g(I+vf)g^{-1}=I+(gv)(f\circ g^{-1}).
\]
The functional $f\circ g^{-1}$ is primitive, and $g$ maps
$\ker f$ bijectively onto $\ker(f\circ g^{-1})$.
Letting $v$ range over $\ker f$ therefore gives
\begin{equation}\label{eq:conjugate-U}
 gU_fg^{-1}=U_{f\circ g^{-1}}.
\end{equation}

\begin{lemma}\label{lem:dual-basis-U}
If $f_1,\ldots,f_r$ form an $R$-basis of $V^*$, then the
subgroups $U_{f_1},\ldots,U_{f_r}$ generate $\SL_r(R)$.
\end{lemma}
\begin{proof}
Let $e_1,\ldots,e_r$ be the dual basis of $V$.
For $i\neq j$ and $a\in R$, the vector $ae_j$ lies in
$\ker f_i$, so $I+ae_jf_i\in U_{f_i}$. These are all
elementary matrices in the basis $e_1,\ldots,e_r$.
Apply Lemma~\ref{lem:elementary-generation}.
\end{proof}

\begin{proposition}\label{prop:elementary-image}
Let $R$ be a commutative local ring with residue field $k$,
and let $G\subseteq\GL_r(R)$. Suppose the reduction of $G$
acts irreducibly on $k^r$ and $G$ contains $U_f$ for a
primitive functional $f$. Then $G\supseteq\SL_r(R)$.
If also $\det G=R^\times$, then $G=\GL_r(R)$.
\end{proposition}
\begin{proof}
The dual residual representation is irreducible. Indeed, the
annihilator of a proper nonzero invariant subspace of $(k^r)^*$
would be a proper nonzero invariant subspace of $k^r$.
Since $f$ is primitive, its reduction $\overline f$ is nonzero.
The span of its orbit is therefore the whole dual space.
Choose $g_1,\ldots,g_r\in G$ such that the reductions of
$f_i=f\circ g_i^{-1}$ form a basis of $(k^r)^*$.

The matrix with rows $f_1,\ldots,f_r$ has nonzero determinant
modulo the maximal ideal. Its determinant is a unit in $R$,
so these functionals form an $R$-basis of $V^*$.
By \eqref{eq:conjugate-U}, $G$ contains all the groups $U_{f_i}$.
Lemma~\ref{lem:dual-basis-U} now implies that $\SL_r(R)\subseteq G$.

If $\det G=R^\times$, take any $M\in\GL_r(R)$ and choose
$g\in G$ with $\det g=\det M$. Then
$Mg^{-1}\in\SL_r(R)\subseteq G$, whence $M\in G$.
\end{proof}

\begin{proposition}\label{prop:local-from-irreducible}
Suppose $\overline\rho_{\phi,\lambda}$ is irreducible and
there is a prime $\mathfrak p\neq\lambda$ with
\[
 v_{\mathfrak p}(g_{r-1})=0\quad\text{and}\quad v_{\mathfrak p}(g_r)=1.
\]
Then $\widehat\rho_{\phi,\lambda}(G_F)=\GL_r(A_\lambda)$.
\end{proposition}
\begin{proof}
Theorem~\ref{thm:full-U} implies that the full translation group
$U_f$ lies inside the image over $A_\lambda$.
Proposition~\ref{prop:elementary-image} therefore asserts
$\SL_r(A_\lambda)$ is contained inside that image.
Corollary~\ref{cor:det-local} gives surjectivity for the determinant
$A_\lambda^\times$, so the last assertion of
Proposition~\ref{prop:elementary-image} proves the claim.
\end{proof}

\section{Non-entanglement of the $\lambda$-adic images}\label{sec:independence}

We now study possible entanglement between the $\lambda$-adic
images $\widehat\rho_{\phi,\lambda}(G_F)\subseteq\GL_r(A_\lambda)$.
For distinct primes $\lambda_1,\lambda_2$, write
\[
 G_{12}
 =(\widehat\rho_{\phi,\lambda_1},
   \widehat\rho_{\phi,\lambda_2})(G_F)\subseteq \GL_r(A_{\lambda_1})\times \GL_r(A_{\lambda_2})
\]
for the joint image. We show that, if both $\lambda$-adic images are full and a prime
$\mathfrak p\notin\{\lambda_1,\lambda_2\}$ satisfies
\eqref{eq:intro-stable}, then
\[
 G_{12}\supseteq
 \SL_r(A_{\lambda_1})\times\SL_r(A_{\lambda_2}).
\]
The proof uses Theorem~\ref{thm:full-U} as well as \cite[Lemma~3.3]{ZywinaAverage}.

\begin{lemma}\label{lem:two-adic-U}
Suppose $\mathfrak p$ satisfies \eqref{eq:intro-stable}, and
let $\lambda_1,\lambda_2$ be distinct primes different from
$\mathfrak p$. In compatible bases adapted to uniformization,
the image of $P_{F_{\mathfrak p}}$ in the two local factors is
\[
 U(A_{\lambda_1})\times U(A_{\lambda_2}).
\]
\end{lemma}
\begin{proof}
For $i=1,2$, the uniformization sequence is
\[
 0\longrightarrow T_{\lambda_i}(\psi)
 \longrightarrow T_{\lambda_i}(\phi)
 \longrightarrow \Lambda\otimes_A A_{\lambda_i}
 \longrightarrow0.
\]
Choose a basis of $T_{\lambda_i}(\psi)$ and extend it to a
basis of $T_{\lambda_i}(\phi)$ by lifting a generator of the
rank-one quotient. Reducing these bases modulo powers of
$\lambda_i$ gives compatible adapted bases at every finite level.

Since $\psi$ has good reduction at $\mathfrak p$ and
$\lambda_i\neq\mathfrak p$, inertia acts trivially on
$T_{\lambda_i}(\psi)$. Moreover, the action on the rank-one
lattice $\Lambda$ takes values in $A^\times=\F_q^\times$.
This group has order prime to $p$, so the pro-$p$ group
$P_{F_{\mathfrak p}}$ acts trivially on $\Lambda$ and hence
on $\Lambda\otimes_A A_{\lambda_i}$. Thus, in the chosen bases,
the joint image
\[
 H=
 (\widehat\rho_{\phi,\lambda_1},
  \widehat\rho_{\phi,\lambda_2})(P_{F_{\mathfrak p}})
\]
satisfies
\[
 H\subseteq U(A_{\lambda_1})\times U(A_{\lambda_2}).
\]

We prove that this inclusion is an equality by considering
both factors simultaneously at finite level. Fix $n,m\geq1$
and put $\mathfrak a=\lambda_1^n\lambda_2^m$.
The Chinese remainder theorem gives a Galois-equivariant
decomposition
\[
 \phi[\mathfrak a]\simeq
 \phi[\lambda_1^n]\oplus\phi[\lambda_2^m],
\]
compatible with the submodules arising from $\psi$.
The two reduced adapted bases therefore determine an adapted
basis of $\phi[\mathfrak a]$ over $A/\mathfrak a$.

Since $\mathfrak a$ is prime to $\mathfrak p$,
Theorem~\ref{thm:full-U} identifies the wild inertia image
at this level with $U(A/\mathfrak a)$. This group is
\[
 U(A/\mathfrak a)
 \simeq U(A/\lambda_1^n)\times U(A/\lambda_2^m).
\]
Consequently, reduction induces a surjection
\[
 H\longrightarrow
 U(A/\lambda_1^n)\times U(A/\lambda_2^m).
\]
In particular, any prescribed pair of translation matrices
modulo $\lambda_1^n$ and $\lambda_2^m$ is realized by the
same element of $P_{F_{\mathfrak p}}$.

As $n$ and $m$ vary, the kernels of these reduction maps
form a neighbourhood basis of the identity in
$U(A_{\lambda_1})\times U(A_{\lambda_2})$. The preceding
surjectivity therefore shows that $H$ is dense in this product.
On the other hand, $P_{F_{\mathfrak p}}$ is compact and the
representations are continuous, so $H$ is compact and hence
closed in the Hausdorff target. It follows that
\[
 H=U(A_{\lambda_1})\times U(A_{\lambda_2}),
\]
as required.
\end{proof}

\begin{proposition}\label{prop:pairwise}
Let $\lambda_1,\lambda_2$ be distinct primes such that
\[
 \widehat\rho_{\phi,\lambda_i}(G_F)=\GL_r(A_{\lambda_i})
 \qquad(i=1,2).
\]
Suppose there is a prime
$\mathfrak p\notin\{\lambda_1,\lambda_2\}$ satisfying
\eqref{eq:intro-stable}. Then the joint image
\[
 G_{12}
 =(\widehat\rho_{\phi,\lambda_1},
   \widehat\rho_{\phi,\lambda_2})(G_F)
\]
contains
$\SL_r(A_{\lambda_1})\times\SL_r(A_{\lambda_2})$.
\end{proposition}

\begin{proof}
Put $R_i=A_{\lambda_i}$ and $S_i=\SL_r(R_i)$.
By Lemma~\ref{lem:two-adic-U}, there are primitive functionals
$f_i:(R_i)^r\to R_i$ such that
\[
 U_{f_1}\times U_{f_2}\subseteq G_{12}.
\]
In particular,
\[
 U_{f_1}\times\{I\}\subseteq G_{12},
 \qquad
 \{I\}\times U_{f_2}\subseteq G_{12}.
\]
We first show that the first of these inclusions implies
$S_1\times\{I\}\subseteq G_{12}$.

Let $\epsilon_1,\ldots,\epsilon_r$ be the coordinate
functionals on $R_1^r$. Since $f_1$ is primitive, there is
a basis in which it is a coordinate functional. Permuting
this basis shows that, for every $j$, there is a matrix
$h_j\in\GL_r(R_1)$ such that
\[
 f_1\circ h_j^{-1}=\epsilon_j.
\]
The first projection of $G_{12}$ is $\GL_r(R_1)$ by
hypothesis. We can therefore choose
\[
 g_j=(h_j,k_j)\in G_{12}
\]
for some $k_j\in\GL_r(R_2)$.

Conjugating by $g_j$ and using \eqref{eq:conjugate-U}, we obtain
\[
 g_j\bigl(U_{f_1}\times\{I\}\bigr)g_j^{-1}
 =U_{\epsilon_j}\times\{I\}
 \subseteq G_{12}.
\]
The second component remains the identity because
$k_jIk_j^{-1}=I$. Thus $G_{12}$ contains all the groups
$U_{\epsilon_j}\times\{I\}$, for $1\leq j\leq r$.
Since the $\epsilon_j$ form a basis of $(R_1^r)^*$,
Lemma~\ref{lem:dual-basis-U} shows that these groups generate
$S_1\times\{I\}$.

Applying the same argument to the second projection gives
$\{I\}\times S_2\subseteq G_{12}$. Finally, for arbitrary
$s_i\in S_i$,
\[
 (s_1,s_2)=(s_1,I)(I,s_2)\in G_{12}.
\]
Hence we conclude that $S_1\times S_2\subseteq G_{12}$.
\end{proof}

\begin{theorem}\label{thm:adelic-independent}
Suppose that
\[
 \widehat\rho_{\phi,\lambda}(G_F)=\GL_r(A_\lambda)
\]
for every nonzero prime ideal $\lambda\subset A$, and that at least
three distinct primes satisfy \eqref{eq:intro-stable}.
Let $G=\widehat\rho_\phi(G_F)$. Then one has that
\[
 [G,G]=\SL_r(\widehat A),
 \qquad
 G=(\det)^{-1}(\det G),
\]
where $[G,G]$ denotes the commutator subgroup of $G$.
\end{theorem}

\begin{proof}
Set
\[
 S_\lambda=\SL_r(A_\lambda),
 \qquad
 S=\prod_\lambda S_\lambda=\SL_r(\widehat A), \quad \text{and}\quad
 H=\overline{[G,G]}.
\]
The closure is taken in $\GL_r(\widehat A)$.
Since $G$ is the continuous image of the compact group $G_F$,
it is compact and hence closed. Therefore $H\subseteq G$.
Every commutator has determinant one, and $S$ is closed, so
$H\subseteq S$.

We first determine the projection of $H$ onto any two factors.
Fix distinct primes $\lambda_1,\lambda_2$, and write
\[
 \pi_{12}:G\longrightarrow G_{12}
\]
for the projection onto the corresponding joint image.
Of the three primes satisfying \eqref{eq:intro-stable},
at least one differs from both $\lambda_1$ and $\lambda_2$.
Proposition~\ref{prop:pairwise} therefore gives
\[
 S_{\lambda_1}\times S_{\lambda_2}\subseteq G_{12}.
\]
Each $S_{\lambda_i}$ is perfect by Lemma~\ref{lem:perfect},
so their product is also perfect. It follows that
\[
 \begin{aligned}
 S_{\lambda_1}\times S_{\lambda_2}
 &=[S_{\lambda_1}\times S_{\lambda_2},
    S_{\lambda_1}\times S_{\lambda_2}]\\
 &\subseteq [G_{12},G_{12}]
 \subseteq S_{\lambda_1}\times S_{\lambda_2}.
 \end{aligned}
\]
Thus
\[
 [G_{12},G_{12}]=S_{\lambda_1}\times S_{\lambda_2}.
\]

Since $\pi_{12}$ is surjective,
\[
 \pi_{12}([G,G])=[G_{12},G_{12}].
\]
Indeed, a commutator $[x,y]$ in $G_{12}$ is the image of
$[\widetilde x,\widetilde y]$ for any lifts
$\widetilde x,\widetilde y\in G$, and the same observation
applies to finite products of commutators. Since
$[G,G]\subseteq H\subseteq S$, we conclude that
\[
 \pi_{12}(H)=S_{\lambda_1}\times S_{\lambda_2}.
\]

Now let $\mathcal L$ be a finite set of at least two primes,
and let $H_{\mathcal L}$ be the image of $H$ in
$\prod_{\lambda\in\mathcal L}S_\lambda$.
The group $H_{\mathcal L}$ is compact and hence closed.
Its projection onto every pair of factors is surjective by
the preceding argument. Each factor is perfect, so
\cite[Lemma~3.3]{ZywinaAverage} gives
\[
 H_{\mathcal L}=\prod_{\lambda\in\mathcal L}S_\lambda.
\]

Consequently, the projection of $H$ onto every finite
collection of factors of $S$ is surjective. To see that this
makes $H$ dense in $S$, take any $s=(s_\lambda)_\lambda\in S$
and any basic open neighbourhood of $s$. Such a neighbourhood
imposes conditions on only finitely many coordinates.
Surjectivity onto those coordinates supplies an element of
$H$ agreeing with $s$ there, and hence lying in the
neighbourhood. Since $H$ is closed, we obtain
\[
 H=S=\SL_r(\widehat A).
\]

In particular, $S\subseteq G$. We now use the perfectness of $S$, established in Lemma~\ref{lem:perfect}:
\[
 S=[S,S]\subseteq[G,G]\subseteq S.
\]
Thus $[G,G]=S$.

It remains to determine $G$ from its determinant image.
The inclusion $G\subseteq(\det)^{-1}(\det G)$ follows from the
definition. Conversely, take
$M\in\GL_r(\widehat A)$ with $\det M\in\det G$.
Choose $g\in G$ such that $\det g=\det M$. Then
\[
 \det(Mg^{-1})=1,
\]
so $Mg^{-1}\in S\subseteq G$, and hence $M\in G$.
This proves $G=(\det)^{-1}(\det G)$.
\end{proof}

\section{Congruence conditions and density estimates}\label{sec:sieve}

\par We begin by recalling that $\B=A^{r-1}\times(A\setminus\{0\})$ and that $\B_n=A_{\leq n}^{r-1}\times(A_{\leq n}\setminus\{0\})$. In particular,
\[
 |\B_n|
 =q^{(r-1)(n+1)}(q^{n+1}-1)
 =q^{r(n+1)}\bigl(1-q^{-(n+1)}\bigr).
\]

Let $\mathfrak b$ be a fixed nonzero ideal of degree $e$.
For $n+1\geq e$, reduction gives a surjective
$\F_q$-linear map $A_{\leq n}\longrightarrow A/\mathfrak b$.
Its source and target have dimensions $n+1$ and $e$,
respectively as $\F_q$-vector spaces. Every fibre therefore has
$q^{n+1-e}$ elements. Consequently, for any subset
$\Omega\subseteq(A/\mathfrak b)^r$, the proportion of tuples
in $A_{\leq n}^r$ reducing into $\Omega$ is exactly
$|\Omega|/|A/\mathfrak b|^r$ once $n+1\geq e$.

The tuples with $g_r=0$ is a set with density $q^{-(n+1)}$ in $A_{\leq n}^r$. Removing them does not change the limiting values. Thus we find that
\[
 \lim_{n\to\infty}
 \frac{|\{\boldsymbol g\in\B_n:
              \boldsymbol g\bmod\mathfrak b\in\Omega\}|}
      {|\B_n|}
 =\frac{|\Omega|}{|A/\mathfrak b|^r}.
\]

\begin{definition}\label{def:stable-prime-condition}
For a nonzero prime $\mathfrak p\subset A$, let
\[
 E_{\mathfrak p}
 =\left\{\boldsymbol g\in\B:
 v_{\mathfrak p}(g_{r-1})=0,\quad
 \text{and}\quad v_{\mathfrak p}(g_r)=1\right\}.
\]

These conditions depend only on the tuple modulo
$\mathfrak p^2$. Let
$\overline E_{\mathfrak p}\subseteq(A/\mathfrak p^2)^r$
be the set of residue tuples satisfying them, and define
\[
 c_{\mathfrak p}
 =\frac{|\overline E_{\mathfrak p}|}
        {|A/\mathfrak p^2|^r}.
\]
In other words, $c_{\mathfrak p}$ is the probability that
the conditions defining $E_{\mathfrak p}$ hold when the
coefficients are chosen uniformly modulo $\mathfrak p^2$.
\end{definition}

\begin{lemma}\label{lem:stable-prime-proportions}
For every nonzero prime $\mathfrak p$, the set
$E_{\mathfrak p}$ has density $c_{\mathfrak p}$. Writing
$Q_{\mathfrak p}=\Nm\mathfrak p=q^{\deg\mathfrak p}$, we have
\begin{equation}\label{eq:c-p}
 c_{\mathfrak p}
 =\left(1-Q_{\mathfrak p}^{-1}\right)
  \left(Q_{\mathfrak p}^{-1}-Q_{\mathfrak p}^{-2}\right)
 =Q_{\mathfrak p}^{-1}
  \left(1-Q_{\mathfrak p}^{-1}\right)^2.
\end{equation}
Moreover,
\[
 \sum_{\mathfrak p}c_{\mathfrak p}=\infty.
\]
There are also at least two primes of every sufficiently
large degree.
\end{lemma}

\begin{proof}
Fix $\mathfrak p$, and put $Q=Q_{\mathfrak p}$.
We first compute the proportion $c_{\mathfrak p}$.
There are $Q^2$ residue classes modulo $\mathfrak p^2$,
of which $Q$ are divisible by $\mathfrak p$. Hence there
are $Q^2-Q$ choices for $g_{r-1}$ satisfying
$v_{\mathfrak p}(g_{r-1})=0$.

For $g_r$, we require divisibility by $\mathfrak p$ but
not by $\mathfrak p^2$. Its residue class must therefore
be a nonzero element of $\mathfrak p/\mathfrak p^2$,
giving $Q-1$ choices. The other $r-2$ coefficients are
unrestricted, with $Q^2$ choices each. Thus
\[
 |\overline E_{\mathfrak p}|
 =(Q^2-Q)(Q-1)Q^{2(r-2)}.
\]
Dividing by the total number $Q^{2r}$ of residue tuples gives
\[
 c_{\mathfrak p}
 =\frac{Q^2-Q}{Q^2}\,\frac{Q-1}{Q^2}
 =Q^{-1}(1-Q^{-1})^2.
\]

For all sufficiently large $n$, the tuples in
$A_{\leq n}^r$ are uniformly distributed modulo
$\mathfrak p^2$. Exactly a proportion $c_{\mathfrak p}$
satisfies the required conditions. Every such tuple has
$g_r\neq0$, so it already belongs to $\B_n$. Consequently,
\[
 \frac{|E_{\mathfrak p}\cap\B_n|}{|\B_n|}
 =\frac{c_{\mathfrak p}q^{r(n+1)}}
        {q^{r(n+1)}(1-q^{-(n+1)})}
 \longrightarrow c_{\mathfrak p}.
\]
This proves $\delta(E_{\mathfrak p})=c_{\mathfrak p}$.

We now prove that $\sum_{\mathfrak p}c_{\mathfrak p}$
diverges. The formula just obtained shows that a prime
of degree $d$ contributes
\[
 q^{-d}(1-q^{-d})^2.
\]
To determine the total contribution from degree $d$, we
must therefore count how many primes have that degree.
Let $\pi_q(d)$ denote this number, or equivalently the
number of monic irreducible polynomials of degree $d$.
Then
\[
 \sum_{\deg\mathfrak p=d}c_{\mathfrak p}
 =\pi_q(d)q^{-d}(1-q^{-d})^2.
\]
We will show that this quantity is asymptotic to $1/d$.

An irreducible polynomial of degree $d$ over $\F_q$
divides $T^{q^n}-T$ precisely when $d\mid n$.
Moreover, $T^{q^n}-T$ has no repeated factors, since its
derivative is $-1$. Comparing degrees in its factorization
therefore gives
\[
 q^n=\sum_{d\mid n}d\pi_q(d).
\]
By M\"obius inversion,
\[
 d\pi_q(d)=\sum_{a\mid d}\mu(a)q^{d/a},
\]
where $\mu$ is the M\"obius function. The term with $a=1$
is $q^d$. For every other divisor $a$, we have $d/a\leq d/2$,
so
\[
 \left|
 \sum_{\substack{a\mid d\\a>1}}\mu(a)q^{d/a}
 \right|
 \leq\sum_{j=1}^{\lfloor d/2\rfloor}q^j
 =O_q(q^{d/2}).
\]
It follows that
\[
 \pi_q(d)=\frac{q^d}{d}
          +O_q\left(\frac{q^{d/2}}d\right).
\]

Substituting this estimate into the contribution from
degree $d$, we obtain
\[
 \sum_{\deg\mathfrak p=d}c_{\mathfrak p}
 =\pi_q(d)q^{-d}(1-q^{-d})^2
 =\frac1d+O_q\left(\frac{q^{-d/2}}d\right).
\]
Thus the number of primes of degree $d$ compensates for
the small contribution of each individual prime.
Summing over $1\leq d\leq D$ gives
\[
 \sum_{\deg\mathfrak p\leq D}c_{\mathfrak p}
 =\sum_{d=1}^{D}\frac1d+O_q(1),
\]
because $\sum_{d\geq1}q^{-d/2}/d$ converges.
The harmonic sum tends to infinity, proving
$\sum_{\mathfrak p}c_{\mathfrak p}=\infty$.

Finally, the estimate for $\pi_q(d)$ shows that
$\pi_q(d)\to\infty$. In particular, there are at least
two primes of every sufficiently large degree.
\end{proof}

We now use the divergence of $\sum_{\mathfrak p}c_{\mathfrak p}$
to show that almost every coefficient tuple satisfies the
required valuation conditions at at least three primes.

\begin{proposition}\label{prop:stable-density}
Let $\cR$ be the set of tuples satisfying $E_{\mathfrak p}$
at at least three distinct primes. Then $\delta(\cR)=1$.
\end{proposition}

\begin{proof}
Let $\varepsilon>0$. Since
$\sum_{\mathfrak p}c_{\mathfrak p}=\infty$, we can choose
a finite set $\mathcal P$ of primes such that
\[
 S_{\mathcal P}:=\sum_{\mathfrak p\in\mathcal P}
                         c_{\mathfrak p}>2\quad \text{and}\quad
 \frac{S_{\mathcal P}}{(S_{\mathcal P}-2)^2}<\varepsilon.
\]
We keep this set $\mathcal P$ fixed throughout the following
counting argument.

Define
\[
 \mathcal C_{\mathcal P}
 =\left\{\boldsymbol g\in\B:
   \#\{\mathfrak p\in\mathcal P:
             \boldsymbol g\in E_{\mathfrak p}\}\leq2
   \right\}.
\]
This set depends on $\mathcal P$, but it contains the fixed
exceptional set $\B\setminus\cR$. Indeed, a tuple satisfying
the valuation conditions at fewer than three primes in total
necessarily satisfies them at fewer than three primes in
$\mathcal P$. Thus
\[
 \B\setminus\cR\subseteq\mathcal C_{\mathcal P}.
\]

We next bound the density of $\mathcal C_{\mathcal P}$.
Put
\[
 \mathfrak b_{\mathcal P}
 =\prod_{\mathfrak p\in\mathcal P}\mathfrak p^2,
\]
and choose a residue tuple uniformly in
$(A/\mathfrak b_{\mathcal P})^r$.
For each $\mathfrak p\in\mathcal P$, let $X_{\mathfrak p}$
be one if the tuple satisfies the conditions defining
$E_{\mathfrak p}$, and zero otherwise. Set
\[
 X_{\mathcal P}
 =\sum_{\mathfrak p\in\mathcal P}X_{\mathfrak p}.
\]
The Chinese remainder theorem gives
\[
 (A/\mathfrak b_{\mathcal P})^r
 \simeq
 \prod_{\mathfrak p\in\mathcal P}(A/\mathfrak p^2)^r.
\]
Under the uniform measure, $X_{\mathfrak p}$ depends only on the
$\mathfrak p$-component and has mean $c_{\mathfrak p}$.
Consequently,
\[
 \mathbb E[X_{\mathcal P}]=S_{\mathcal P}\quad \text{and}\quad 
 \operatorname{Var}(X_{\mathcal P})
 =\sum_{\mathfrak p\in\mathcal P}
       c_{\mathfrak p}(1-c_{\mathfrak p})
 \leq S_{\mathcal P}.
\]
If $X_{\mathcal P}\leq2$, then
$|X_{\mathcal P}-S_{\mathcal P}|\geq S_{\mathcal P}-2$.
Chebyshev's inequality therefore gives
\[
 \Pr(X_{\mathcal P}\leq2)
 \leq
 \frac{\operatorname{Var}(X_{\mathcal P})}
      {(S_{\mathcal P}-2)^2}
 \leq
 \frac{S_{\mathcal P}}{(S_{\mathcal P}-2)^2}
 <\varepsilon.
\]

Membership in $\mathcal C_{\mathcal P}$ is determined
modulo the fixed ideal $\mathfrak b_{\mathcal P}$.
For $n+1\geq\deg\mathfrak b_{\mathcal P}$, reduction from
$A_{\leq n}^r$ to $(A/\mathfrak b_{\mathcal P})^r$ has
equal-sized fibres. The proportion of tuples in
$A_{\leq n}^r$ satisfying at most two of the conditions
is therefore exactly $\Pr(X_{\mathcal P}\leq2)$.
Removing the tuples with $g_r=0$ does not change the limit,
since their proportion is $q^{-(n+1)}$. Hence
\[
 \delta(\mathcal C_{\mathcal P})
 =
 \lim_{n\to\infty}
 \frac{|\mathcal C_{\mathcal P}\cap\B_n|}{|\B_n|}
 =
 \Pr(X_{\mathcal P}\leq2)
 <\varepsilon.
\]

Using $\B\setminus\cR\subseteq\mathcal C_{\mathcal P}$,
we obtain
\[
 \overline\delta(\B\setminus\cR)
 \leq\delta(\mathcal C_{\mathcal P})
 <\varepsilon.
\]
The set $\B\setminus\cR$ is independent of $\varepsilon$
and $\mathcal P$. Since $\varepsilon>0$ was arbitrary,
its upper density is zero. Therefore $\delta(\cR)=1$.
\end{proof}

Let $\cS_m$ be the set of tuples for which no prime of degree
greater than $m$ divides both $g_{r-1}$ and $g_r$. We set $\varepsilon_m
 =\overline\delta(\B\setminus\cS_m)$.

\begin{lemma}\label{lem:common-tail}We have that $\varepsilon_m\to0$ as $m\to\infty$.
\end{lemma}

\begin{proof}
Let $\mathcal D_m\subseteq A^2$ consist of pairs $(u,v)$
such that some prime of degree greater than $m$ divides
both $u$ and $v$. By \cite[Lemma~9.5]{ZywinaAverage},
\[
 \eta_m:=
 \limsup_{n\to\infty}
 \frac{|\mathcal D_m\cap A_{\leq n}^2|}{q^{2(n+1)}}
 \longrightarrow0
\]
as $m\to\infty$. A tuple in $\B_n\setminus\cS_m$ has
$(g_{r-1},g_r)\in\mathcal D_m\cap A_{\leq n}^2$.
For each such pair, there are $q^{(r-2)(n+1)}$ choices
for the other coordinates. Consequently,
\[
 |\B_n\setminus\cS_m|
 \leq q^{(r-2)(n+1)}
       |\mathcal D_m\cap A_{\leq n}^2|.
\]
Dividing by $|\B_n|$ gives
\[
 \frac{|\B_n\setminus\cS_m|}{|\B_n|}
 \leq
 \frac{1}{1-q^{-(n+1)}}
 \frac{|\mathcal D_m\cap A_{\leq n}^2|}{q^{2(n+1)}}.
\]
Taking the limit superior in $n$, with $m$ fixed, yields
$\varepsilon_m\leq\eta_m$. This completes the proof.
\end{proof}

For all sufficiently large integers $m$, set
\[
 d(m)=\left\lfloor\frac{m}{2C_r}\right\rfloor\quad\text{and}\quad
 \mathfrak a_m=\prod_{\deg\lambda\leq m}\lambda.
\]
Since $d(m)\to\infty$, we may choose two distinct primes
$\mathfrak p_1(m),\mathfrak p_2(m)$ of degree $d(m)$.
These primes are fixed once $m$ is chosen and do not depend
on the coefficient tuple.

Let $\cT_m$ be the set of tuples for which
\[
 \mathfrak p_1(m)\nmid g_r\quad\text{and}\quad 
 \mathfrak p_2(m)\nmid g_r.
\]
For a fixed prime $\mathfrak p$ of degree $d$ and $n\geq d$,
there are $q^{n+1-d}-1$ nonzero multiples of its monic
generator in $A_{\leq n}$. Thus the proportion of nonzero
polynomials in $A_{\leq n}$ divisible by $\mathfrak p$ is
\[
 \frac{q^{n+1-d}-1}{q^{n+1}-1}
 \longrightarrow q^{-d}.
\]
The other coefficients can be chosen freely, so the same
limiting proportion holds in $\B_n$. Applying the union
bound to the two chosen primes gives
\begin{equation}\label{eq:good-tail}
 \overline\delta(\B\setminus\cT_m)
 \leq2q^{-d(m)}.
\end{equation}
Every tuple in $\cT_m$ has good reduction at both primes:
all coefficients are integral there, and the leading
coefficient $g_r$ is a unit.

Let
\[
 \cU_m=
 \left\{\boldsymbol g\in\B:
 \rho_{\boldsymbol g,\mathfrak a_m}(G_F)
       =\GL_r(A/\mathfrak a_m)\right\}.
\]
Theorem~\ref{cor:fixed-level} implies that
$\delta(\cU_m)=1$. Define
\begin{equation}\label{eq:V-m}
 \cV_m=\cR\cap\cS_m\cap\cT_m\cap\cU_m.
\end{equation}

\begin{proposition}\label{prop:V-density}
For all sufficiently large $m$,
\[
 \underline\delta(\cV_m)
 \geq1-\varepsilon_m-2q^{-d(m)}.
\]
Consequently,
$\underline\delta(\cV_m)\to1$ as $m\to\infty$.
\end{proposition}

\begin{proof}
Given a fixed value $m$, a tuple outside $\cV_m$ fails at least one of
the four conditions in \eqref{eq:V-m}. Therefore
\[
 \B\setminus\cV_m
 \subseteq
 (\B\setminus\cR)\cup
 (\B\setminus\cS_m)\cup
 (\B\setminus\cT_m)\cup
 (\B\setminus\cU_m).
\]
Counting in $\B_n$, dividing by $|\B_n|$, and taking
the limit superior gives
\[
 \begin{aligned}
 \overline\delta(\B\setminus\cV_m)
 &\leq
 \overline\delta(\B\setminus\cR)
 +\overline\delta(\B\setminus\cS_m)\\
 &\quad+
 \overline\delta(\B\setminus\cT_m)
 +\overline\delta(\B\setminus\cU_m).
 \end{aligned}
\]
The first and fourth terms are zero. The second is
$\varepsilon_m$, and the third is at most $2q^{-d(m)}$.
Since
\[
 \underline\delta(\cV_m)
 =1-\overline\delta(\B\setminus\cV_m),
\]
we obtain the asserted lower bound. Finally,
$\varepsilon_m\to0$ by Lemma~\ref{lem:common-tail},
and $q^{-d(m)}\to0$ because $d(m)\to\infty$.
\end{proof}

\section{Proof of the main theorems}\label{sec:main}

\par In this section, we prove the main results of the article.
\begin{proposition}\label{prop:V-maximal}
For every sufficiently large $m$ and every
$\boldsymbol g\in\cV_m$,
\[
 \widehat\rho_{\boldsymbol g,\lambda}(G_F)
 =\GL_r(A_\lambda)
\]
for every nonzero prime $\lambda\subset A$.
\end{proposition}

\begin{proof}
Fix $\boldsymbol g\in\cV_m$, and let $\phi$ be the
Drinfeld module determined by this tuple. We first prove
that $\overline\rho_{\boldsymbol g,\lambda}$ is irreducible
for every $\lambda$.

Suppose first that $\deg\lambda\leq m$. Then
$\lambda\mid\mathfrak a_m$. Since $\boldsymbol g\in\cU_m$,
the representation at level $\mathfrak a_m$ has image
$\GL_r(A/\mathfrak a_m)$. The Chinese remainder theorem
shows that the projection
\[
 \GL_r(A/\mathfrak a_m)\longrightarrow\GL_r(k_\lambda)
\]
is surjective. Hence
\[
 \overline\rho_{\boldsymbol g,\lambda}(G_F)
 =\GL_r(k_\lambda).
\]
The natural action of this group on $k_\lambda^r$ is
irreducible: the orbit of any nonzero vector consists of
all nonzero vectors, so a nonzero invariant subspace must
be the entire space.

Now suppose that $\deg\lambda>m$.
Since $\boldsymbol g\in\cS_m$, the prime $\lambda$ does not
divide both $g_{r-1}$ and $g_r$. There are two possibilities.
If $\lambda\nmid g_r$, all coefficients are integral at
$\lambda$ and the leading coefficient is a unit, so $\phi$
has good reduction there. If $\lambda\mid g_r$, then
$\lambda\nmid g_{r-1}$. The reduction of the displayed
integral model then has rank $r-1$, so $\phi$ has stable
reduction of that rank at $\lambda$.

Membership in $\cT_m$ supplies the two distinct good primes
$\mathfrak p_1(m),\mathfrak p_2(m)$, both of degree $d(m)$.
Thus the reduction hypotheses of
Theorem~\ref{thm:exceptional-bound} hold. If
$\overline\rho_{\boldsymbol g,\lambda}$ were reducible,
that theorem would give
\[
 \deg\lambda
 \leq C_rd(m)
 =C_r\left\lfloor\frac{m}{2C_r}\right\rfloor
 \leq\frac m2.
\]
This contradicts $\deg\lambda>m$. The residual
representation is therefore irreducible in this case
as well.

We have proved residual irreducibility at every prime.
Furthermore, $\boldsymbol g\in\cR$ supplies three distinct
primes satisfying \eqref{eq:intro-stable}. For any given
$\lambda$, choose one of them, say $\mathfrak p$, different
from $\lambda$. Then
\[
 v_{\mathfrak p}(g_{r-1})=0,
 \qquad
 v_{\mathfrak p}(g_r)=1.
\]
The hypotheses of
Proposition~\ref{prop:local-from-irreducible} are now
satisfied at $\lambda$. It follows that
\[
 \widehat\rho_{\boldsymbol g,\lambda}(G_F)
 =\GL_r(A_\lambda).
\]
Since $\lambda$ was arbitrary, this holds for every prime.
\end{proof}

\begin{proof}[Proof of Theorem \ref{thm a}]
For $\boldsymbol g\in\B$, write
\[
 G_{\boldsymbol g}
 =\widehat\rho_{\boldsymbol g}(G_F).
\]
Let $\cM$ be the set of tuples for which
\[
 \widehat\rho_{\boldsymbol g,\lambda}(G_F)
 =\GL_r(A_\lambda)
 \quad\text{for every prime }\lambda,
 \qquad
 [G_{\boldsymbol g},G_{\boldsymbol g}]
 =\SL_r(\widehat A).
\]

We first prove that $\delta(\cM)=1$. Fix a sufficiently
large $m$ and take $\boldsymbol g\in\cV_m$.
Proposition~\ref{prop:V-maximal} gives surjectivity of every
$\lambda$-adic representation. Since
$\boldsymbol g\in\cR$, there are also three distinct primes
satisfying \eqref{eq:intro-stable}.
Theorem~\ref{thm:adelic-independent} therefore gives
\[
 [G_{\boldsymbol g},G_{\boldsymbol g}]
 =\SL_r(\widehat A).
\]
Thus $\cV_m\subseteq\cM$ for every sufficiently large $m$.
By Proposition~\ref{prop:V-density},
\[
 \underline\delta(\cM)
 \geq\underline\delta(\cV_m)
 \geq1-\varepsilon_m-2q^{-d(m)}.
\]
Letting $m\to\infty$ gives $\underline\delta(\cM)\geq1$.
Since all box proportions are at most one, this proves
$\delta(\cM)=1$.

Here the limits defining the densities are taken with
$m$ fixed. The subsequent limit in $m$ improves the lower
bound for the single set $\cM$. This is why the fixed-level
density statements suffice to obtain simultaneous
surjectivity at all primes.

Now fix $\boldsymbol g\in\cM$. Since the commutator subgroup
is contained in the group itself,
\[
 \SL_r(\widehat A)\subseteq G_{\boldsymbol g}.
\]
If $M\in\GL_r(\widehat A)$ has determinant in
$\det G_{\boldsymbol g}$, choose $g\in G_{\boldsymbol g}$
with $\det g=\det M$. Then
$Mg^{-1}\in\SL_r(\widehat A)\subseteq G_{\boldsymbol g}$,
so $M\in G_{\boldsymbol g}$. Hence
\[
 G_{\boldsymbol g}
 =(\det)^{-1}(\det G_{\boldsymbol g}).
\]
The determinant map on $\GL_r(\widehat A)$ is surjective,
as is seen from the matrices $\operatorname{diag}(u,1,\ldots,1)$.
It therefore induces an isomorphism
\[
 \GL_r(\widehat A)/G_{\boldsymbol g}
 \simeq
 \widehat A^\times/\det G_{\boldsymbol g}.
\]
In particular, the adelic index is the determinant index.

The determinant representation is the adelic representation
attached to the rank-one Drinfeld module
\[
 \psi_T=T+\Delta\tau,\quad \text{where}\quad
 \Delta=(-1)^{r-1}g_r.
\]
Put $d=\deg g_r=\deg\Delta$. Then
\[
 (-1)^d\lc(\Delta)
 =(-1)^{d+r-1}\lc(g_r),
\]
whose order in $\F_q^\times$ is $e_r(g_r)$.
Applying \eqref{eq:gekeler-adelic} gives
\[
 \begin{aligned}
 [\GL_r(\widehat A):G_{\boldsymbol g}]
 &=[\widehat A^\times:\det G_{\boldsymbol g}]\\
 &=\gcd\left(d-1,\frac{q-1}{e_r(g_r)}\right).
 \end{aligned}
\]
This index divides $q-1$, and all the asserted properties
hold on the density-one set $\cM$.
\end{proof}

\begin{proof}[Proof of Theorem \ref{thm b}]
Recall that $\cP$ consists of the tuples for which
\[
 (-1)^{\deg g_r+r-1}\lc(g_r)
\]
generates $\F_q^\times$. For such a tuple,
$e_r(g_r)=q-1$. The determinant formula
\eqref{eq:gekeler-adelic} therefore gives
\[
 [\widehat A^\times:\det G_{\boldsymbol g}]
 =\gcd(\deg g_r-1,1)=1.
\]
Thus $\det G_{\boldsymbol g}=\widehat A^\times$ for every
$\boldsymbol g\in\cP$.

If $\boldsymbol g$ also belongs to $\cM$, Theorem A gives
\[
 G_{\boldsymbol g}
 =(\det)^{-1}(\det G_{\boldsymbol g})
 =\GL_r(\widehat A).
\]
Consequently,
\[
 \cM\cap\cP\subseteq\cS_{\mathrm{ad}},
\]
where $\cS_{\mathrm{ad}}$ is the set of tuples with
surjective adelic representation.

Put
\[
 c=\frac{\varphi(q-1)}{q-1}>0.
\]
By Lemma~\ref{lem:primitive-density},
\[
 |\cP\cap\B_n|=c|\B_n|
\]
for every $n$. A tuple in $\cP$ whose adelic representation
is not surjective cannot belong to $\cM$. Hence
\[
 (\cP\setminus\cS_{\mathrm{ad}})\cap\B_n
 \subseteq\B_n\setminus\cM.
\]
Dividing by $|\cP\cap\B_n|$ yields
\[
 0\leq
 1-\frac{|\cS_{\mathrm{ad}}\cap\cP\cap\B_n|}
          {|\cP\cap\B_n|}\leq
 \frac{|\B_n\setminus\cM|}{c|\B_n|}.
\]
On the other hand, 
\[ \frac{|\B_n\setminus\cM|}{c|\B_n|}
 \longrightarrow0\]
since $\delta(\cM)=1$. This proves that adelic
surjectivity has relative density one inside $\cP$.

Note that
\[
 \begin{aligned}
 |\cS_{\mathrm{ad}}\cap\B_n|
 &\geq|\cM\cap\cP\cap\B_n|\\
 &\geq|\cP\cap\B_n|-|\B_n\setminus\cM|.
 \end{aligned}
\]
After dividing by $|\B_n|$ and taking the limit inferior,
we obtain
\[
 \underline\delta(\cS_{\mathrm{ad}})
 \geq c
 =\frac{\varphi(q-1)}{q-1}.
\]
This proves the claimed lower bound.
\end{proof}
\bibliographystyle{alpha}
\bibliography{references}
\end{document}